\documentclass[dvips,11pt]{article}
\usepackage{amsmath,amssymb,amsfonts,amscd,amsthm,pstricks,pst-node,bm,xspace}
\def\simarrow{\mathrel{\raise -0.5mm\hbox{$\sim$}}\hspace{-1.8mm}{\rightarrow} } 

\def\bsimarrow{\leftarrow\hspace{-0.7mm}\mathrel{\raise -0.5mm\hbox{$\backsim$}} }

\def\bt{\begin{tabular}}
\def\te{\end{tabular}}

\def\lettrine#1#2#3{\noindent\hangindent#1\hangafter-#2
\hskip-#1\smash{\hbox to #1{#3\hfill}}\ignorespaces}

\newcommand{\To}[1]{\mathop{\to}\limits_{#1}}

\def\BM{\begin{pmatrix}}
\def\EM{\end{pmatrix}}

\def\d=f{\buildrel\hbox{\scriptsize d\'{e}f}\over \Longleftrightarrow}

\def\cit{\text{\it I\hskip -6ptC\/}}

\def\rit{\text{\it I\hskip -2pt  R}}

\def\rl {\rit^{\hskip 1pt\ell}}

\def\Bd{{\text B}}

\def\Ed{{\text E}}

\def\be{\begin{equation}}
\def\ee{\end{equation}}
\def\beqn{\begin{eqnarray}}
\def\eeqn{\end{eqnarray}}
\def\nobeqn{\begin{eqnarray*}}
\def\noeeqn{\end{eqnarray*}}
\def\ba{\left(\begin{array}}
\def\ea{\end{array} \right) }

\def\u{\underline}
\def\o{\overline}
\def\and{\; \mbox{and} \;}

\def\hfl#1#2{\smash{\mathop{\hbox to 12mm{\rightarrowfill}}
\limits^{\scriptstyle #1}_{\scriptstyle #2}}}

\def\cov{\mathop{\rm cov}\nolimits}

\def\mod{\mathop{\rm mod}\nolimits}

\def\Be{\begin{enumerate}}
\def\Ee{\end{enumerate}}

\def\Bena{\begin{enumerate}
\def\labelenumi{\theenumi)}
\def\theenumi{\arabic{enumi}}
\def\labelenumii{\theenumii)}
\def\theenumii{\alph{enumii}}}

\def\Bean{\begin{enumerate}
\def\labelenumii{\theenumii)}
\def\theenumii{\arabic{enumii}}
\def\labelenumi{\theenumi)}
\def\theenumi{\alph{enumi}}}

\def\Bero{\begin{enumerate}
\def\labelenumii{\theenumii)}
\def\theenumii{\arabic{enumii}}
\def\labelenumi{(\theenumi)}
\def\theenumi{\roman{enumi}}}

\def\BeRo{\begin{enumerate}
\def\labelenumii{\theenumii)}
\def\theenumii{\arabic{enumii}}
\def\labelenumi{(\theenumi)}
\def\theenumi{\Roman{enumi}}}

\def\Bi{\vskip 11pt\begin{itemize}\itemsep=18pt}
\def\Ei{\end{itemize}\vskip 11pt}

\def\Bd{\begin{description}}
\def\Ed{\end{description}}

\def\R{\right}
\def\L{\left}

\usepackage[latin1]{inputenc}

\def\wt{\widetilde}

\def\ST{{\rm ST}}
\def\MG{{\rm MG}}
\def\TC{{\rm TC}}
\def\M{{\rm M}}

\def\DMS{{\rm DMS}}
\def\VMS{{\rm VMS}}
\def\LGC{{\rm LGC}}
\def\LGGC{{\rm LGGC}}
\def\int{{\rm int}}
\def\sSOT{{\text{\tiny SOT}}}

\def\Times{{\mathrel{\u\times}}}

\def\prod{\mathop{\Pi}\limits}
\def\sum{\mathop{\Sigma}\limits}

\def\bbf{\boldmath\bfseries}
\def\o{\overline}

\def\Bi{\begin{itemize}}
\def\Ei{\end{itemize}}

\newcommand{\Aa}{\mathbb{A}\,}

\newcommand{\RR}{\mathbb{R}\,}
\renewcommand{\rit}{\RR}

\def\Gal{\operatorname{Gal}}

\def\mod{\operatorname{mod}}

\def\Repsp{\operatorname{Repsp}}
\def\FRepsp{\operatorname{FRepsp}}

\def\GL{\operatorname{GL}}

\def\cov{\operatorname{cov}}
\def\SOT{\operatorname{SOT}}

\def\bM{\begin{matrix}}
\def\eM{\end{matrix}}

\def\lr{left (resp. right) }
\def\rl{right (resp. left) }

\def\resp#1{(resp. #1)}
\def\rresp#1{\qquad \mbox{(resp.} \quad #1\ )}

\def\To{\begin{CD} @>>>\end{CD}}

\def\RL{_{R\times L}}

\def\wt{\widetilde}

\begin{document}

\setcounter{page}{0}
{\pagestyle{empty}
\null\vfill
\begin{center}
{\LARGE $GL(2)$-structures of the Langlands global program}
\vfill
{Christian {\sc Pierre\/}}
\vfill



\null\vfill
\vfill

\begin{abstract}
All kinds of global correspondences of Langlands are evaluated from the functional representation spaces of the algebraic bilinear semigroups $GL_2(L_{\o v}\times L_v)$ with entries in products, right by left, of sets of archimedean increasing completions.

Degenerate singularities on these functional representation spaces can give rise, by versal deformations and blowups of these, to one or two new covering functional representation spaces of $GL_2(\cdot\times\cdot)$ according to the type of considered singularities.

The discovered correspondences of Langlands are associated with singular and nonsingular universal $GL(2)$-structures.

\end{abstract}
\end{center}
\vfill
\eject

\vfill\eject}
\setcounter{page}{1}
\def\thepage{\arabic{page}}
\parindent=0pt 
\section{Introduction}
The generalization of the concept of Dirichlet character led Langlands \cite{Lan1} to formulate the (global) reciprocity conjecture \cite{Lan2} which asserts that:

\begin{quote}
For any irreducible representation $\sigma $ of $\Gal(\o F/F)$ in $\GL_n(\cit)$, there exists a cuspidal (unramified) automorphic representation $\Pi $ of $\GL_n(\Aa_F)$ in such a way that the Artin $L$-function of $\sigma $ agrees with the Langlands $L$-function $\Pi $ at almost every place of the adele ring of a number field $F$ of characteristic $0$ \cite{Kna},  \cite{Shi}.
\end{quote}
\vskip 11pt

Recent advances were realized in the so-called Langlands program essentially by M. Harris and R. Taylor \cite{H-T} in the case of $p$-adic number fields and by L. Lafforgue \cite{Laf} in the case of function fields.
\vskip 11pt

In the frame of the Langlands global program developed by the author \cite{Pie1}, an interesting question consists in {\bbf evaluating all kinds of correspondences of Langlands from the functional representation spaces of $\GL_2(F)$ being affected or not by (degenerate) singularities\/}.

The discovered correspondences are then associated with singular and nonsingular universal $\GL(2)$-structures.  

The basic mathematical frame (considered in chapter 2 of this work) is essentially bilinear in such a way that the {\bbf Langlands global correspondences consist in bijections between the equivalence classes of  Galois representations given by bilinear algebraic semigroups $\GL_2(L_{\o v}\times L_v)$ and the corresponding conjugacy classes of their cuspidal representations\/} where $L_v$ \resp{$L_{\o v}$} denotes the set of \lr increasing archimedean pseudo-ramified completions.  These completions, resulting from the corresponding compactified Galois extensions, are characterized by increasing extension degrees which are integers modulo $N$.

Let $GL_2(L_{\o v}\times L_v)$ be the bilinear algebraic semigroup with entries in the products, right by left, of symmetric completions. It decomposes according to:
\[\GL_2(L_{\o v}\times L_v)=T_2^t(L_{\o v})\times T_2(L_v)\]
where $T_2^t(L_{\o v})$ 
\resp{$T_2(L_{v})$} is a \rl linear semigroup of lower \resp{upper} triangular matrices referring to the lower \resp{upper} half space.

{\bbf The representation space of the algebraic bilinear semigroup $\GL_2(L_{\o v}\times L_v)$ is given by the  $\GL_2(L_{\o v}\times L_v)$-bisemimodule $M_R(L_{\o v})\times M_L(L_v)$ and corresponds to the representation space
$\sigma (W _{\wt L_{\o v}}\times W _{\wt L_{v}})$ of the product, right by left, of global Weil groups
$W _{\wt L_{\o v}}$ and $W _{\wt L_{v}}$\/}, respectively over the sets of extensions
$\wt L_{\o v}$ and $\wt L_{v}$.

The set of differentiable bifunctions
$\{\phi _R(M_{\o v_{\mu ,m_\mu }})\otimes\phi _L(M_{v_{\mu ,m_\mu }})\}_{\mu ,m_\mu }$, $1\le\mu \le\infty $, is the set of bisections of {\bbf the bisemisheaf of rings
$\phi _R(M_R(L_{\o v}))\otimes\phi _L(M_L(L_{v}))$ on the algebraic bilinear semigroup
$M_R(L_{\o v})\otimes M_L(L_{v})$\/}.

This leads to {\bbf the Langlands global correspondence}
\[ \LGC: \quad \sigma (W_{\wt L_{\o v}}\times W_{\wt L_v})\To \Pi (\GL_2(L_{\o v}\times L_v))\]
between the set $\sigma (W_{\wt L_{\o v}}\times W_{\wt L_v})$ of two-dimensional representative subspaces of the Weil global subgroups, given by the algebraic bilinear semigroups
$(M_R(L_{\o v})\otimes M_L(L_v))$, and their cuspidal representation given by
$\Pi (\GL_2(L_{\o v}\times L_v))$ and obtained from a toroidal compactification of the bisemisheaf 
$\phi _R(M_R(L_{\o v}))\otimes\phi _L(M_L(L_{v}))$.

As {\bbf the bisemisheaf $\phi _R(M_R(L_{\o v})\otimes _{L_{\o v}\times L_v}
\phi _L(M_L(L_{v}))$\/} is a  vector bisemispace giving rise to an inner product bisemispace, it {\bbf can generate an orthogonal complement bisemisheaf
$\phi _R^\perp(M_R^\perp(L_{\o v}))\otimes
\phi _L^\perp(M_L^\perp(L_{v}))$ by means of endomorphisms based on Galois antiautomorphisms\/}.

By this way, the following {\bbf Langlands general global correspondence can be stated\/}:
\[\LGGC_\ST: \quad
\sigma (W_{\wt L_{\o v}}\times W_{\wt L_v})\To 
\Pi (\GL_2^{(r)}(L_{\o v}\times L_v))
+\Pi (\GL_2^\perp(L_{\o v}\times L_v))\]
where:\Bi
\item $\Pi (\GL_2^{(r)}(L_{\o v}\times L_v))$ is the cuspidal representation on the reduced algebraic bilinear semigroup $\GL_2^{(r)}(L_{\o v}\times L_v)$;

\item $\Pi (\GL_2^\perp(L_{\o v}\times L_v))$ is the cuspidal representation on the orthogonal complement algebraic bilinear semigroup
$\GL_2^\perp(L_{\o v}\times L_v)$.
\Ei

Consequently, the functional representation  bisemispace \cite{F-H} of 
$\sigma (W_{\wt L_{\o v}}\times W_{\wt L_v})$, given by {\bbf the direct sum of the reduced bisemisheaf\/}
$\phi _R(M_R^{(r)}(L_{\o v})\otimes \phi _L(M^{(r)}_L(L_v))$ over
$\GL_2^{(r)}(L_{\o v}\times L_v)$ (and written in condensed notation
{\bbf $(\wt M^T_{\ST_R}\otimes\wt M^T_{\ST_L})$) and of the orthogonal bisemisheaf\/}
$\phi _R^\perp(M_R^\perp(L_{\o v})\otimes \phi _L^\perp(M_L((L_v))$ over
$\GL_2^\perp(L_{\o v}\times L_v)$ (and written in condensed notation
{\bbf $(\wt M^S_{\ST_R}\otimes\wt M^S_{\ST_L})$), constitutes a nonsingular universal mathematical structure\/}.  It is also a nonsingular universal physical structure corresponding to the space-time string fields of an elementary particle.
\vskip 11pt

{\bbf Degenerate singularities on the time or space shifted bisemisheaf
$(\wt M^{T_p}_{\ST_R}\otimes\wt M^{T_p}_{\ST_L})$ or
$(\wt M^{S_p}_{\ST_R}\otimes\wt M^{S_p}_{\ST_L})$
can give rise, by versal deformations \cite{G-K} and blowups of these, to one or two new covering (shifted) bisemisheaves\/} according to the kind of considered degenerate singularities: this constitutes the content of chapter~3.

More concretely, {\bbf according to the existence or the absence of degenerate singularities on the basic space-time bisemisheaf\/}
$(\wt M^{T_p-S_p}_{\ST_R}\otimes\wt M^{T_p-S_p}_{\ST_L})$, where
$\wt M^{T_p-S_p}_{\ST_R}=\wt M^{T_p}_{\ST_R}+\wt M^{S_p}_{\ST_R}$,
{\bbf three sets of Langlands general global correspondences ($\LGGC$) can be built\/}:
\Bean
\item If there are no degenerate singularities on
$(\wt M^{T_p-S_p}_{\ST_R}\otimes\wt M^{T_p-S_p}_{\ST_L})$, we get the {\bbf one-level $\LGGC_\ST$ correspondences\/} mentioned above.

\item If there are degenerate singularities of corank 1 and codimension inferior to 3 or degenerate singularities of corank 2 and codimension inferior or equal to 3 on
$(\wt M^{T_p-S_p}_{\ST_R}\otimes\wt M^{T_p-S_p}_{\ST_L})$, we get, after a process of desingularizations \cite{Hir}, \cite{Hau},  and toroidal compactification of 
$(\wt M^{T_p-S_p}_{\ST_R}\otimes\wt M^{T_p-S_p}_{\ST_L})$, the {\bbf two-level $\LGGC_{\ST-\MG}$ correspondences\/}:
\begin{multline*}
 \LGGC_{\ST-\MG}: \quad
\sigma (W_{\wt L_{\o v}}\times W_{\wt L_v})+
\sigma (W_{\wt L_{\o v_{S\circ T(1)}}}\times W_{\wt L_{v_{S\circ  T(1)}}})\hspace{4.5cm}\\ \To 
\Pi (\GL_2^{(r)}(L_{\o v}\times L_v))+\Pi (\GL_2^\perp(L_{\o v}\times L_v))\\
 +
\Pi (\GL_2^{(r)}(L_{\o v_{\cov(1)}}\times L_{v_{\cov(1)}}))
+\Pi (\GL_2^\perp(L_{\o v_{\cov(1)}}\times L_{v_{\cov(1)}}))\end{multline*}
where:
\Bi
\item $\sigma (W_{\wt L_{\o v_{S\circ T(1)}}}\times W_{\wt L_{v_{S\circ  T(1)}}})$ corresponds to the first level bisemisheaf
$(\wt M^{T_p-S_p}_{\MG_R}\otimes\wt M^{T_p-S_p}_{\MG_L})$ covering the basic bisemisheaf
$(\wt M^{T_p-S_p}_{\ST_R}\otimes\wt M^{T_p-S_p}_{\ST_L})$;

\item $\Pi (\GL_2^{(r)}(L_{\o v_{\cov(1)}}\times L_{v_{\cov(1)}}))+\Pi (\GL_2^\perp(L_{\o v_{\cov(1)}}\times L_{v_{\cov(1)}}))$ is the cuspidal  representation of the first level covering algebraic bilinear semigroups over
$L_{\o v_{\cov(1)}}\times L_{v_{\cov(1)}}\approx 
L_{\o v_{S\circ T(1)}}\times L_{v_{S\circ T(1)}}$.
\Ei

\item If there are degenerate singularities of corank 1 and codimension 3
or $(\wt M^{T_p-S_p}_{\ST_R}\otimes\wt M^{T_p-S_p}_{\ST_L})$, we get, after a process of desingularizations and toroidal compactifications of the singular bisemisheaves
$(\wt M^{T_p-S_p}_{\ST_R}\otimes\wt M^{T_p-S_p}_{\ST_L})$ and
$(\wt M^{T_p-S_p}_{\MG_R}\otimes\wt M^{T_p-S_p}_{\MG_L})$, {\bbf the three-level\linebreak
$\LGGC_{\ST-\MG-\M}$ correspondences\/}:
\begin{multline*}
\sigma (W_{\wt L_{\o v}}\times W_{\wt L_v})+
\sigma (W_{\wt L_{\o v_{S\circ T(1)}}}\times W_{\wt L_{v_{S\circ  T(1)}}})+
\sigma (W_{\wt L_{\o v_{S\circ T(2)}}}\times W_{\wt L_{v_{S\circ  T(2)}}})\\
\To 
\Pi (\GL_2^{(r)}(L_{\o v}\times L_v))+\Pi (\GL_2^\perp(L_{\o v}\times L_v))\hspace{6cm}\\
+
\Pi (\GL_2^{(r)}(L_{\o v_{\cov(1)}}\times L_{v_{\cov(1)}}))+\Pi (\GL_2^\perp(L_{\o v_{\cov(1)}}\times L_{v_{\cov(1)}}))\\ +
\Pi (\GL_2^{(r)}(L_{\o v_{\cov(2)}}\times L_{v_{\cov(2)}}))
+\Pi (\GL_2^\perp(L_{\o v_{\cov(2)}}\times L_{v_{\cov(2)}}))\end{multline*}
where:\Bi
\item $\sigma (W_{\wt L_{\o v_{S\circ T(2)}}}\times W_{\wt L_{v_{S\circ  T(2)}}})$ corresponds to the second level bisemisheaf
$(\wt M^{T_p-S_p}_{\M_R}\otimes\wt M^{T_p-S_p}_{\M_L})$ covering the first level bisemisheaf
$(\wt M^{T_p-S_p}_{\MG_R}\otimes\wt M^{T_p-S_p}_{\MG_L})$;

\item $\Pi (\GL_2^{(r)}(L_{\o v_{\cov(2)}}\times L_{v_{\cov(2)}}))+\Pi (\GL_2^\perp(L_{\o v_{\cov(2)}}\times L_{v_{\cov(2)}}))$ is the cuspidal representation of the second level covering algebraic bilinear semigroups.
\Ei
\Ee

It then results that:
\Bena
\item {\bbf the one-level correspondences $\LGGC_{\ST}$ are related to nonsingular universal mathematical and physical $\GL(2)$-structures\/} given by the bisemisheaves
$(\wt M^{T_p-S_p}_{\ST_R}\otimes\wt M^{T_p-S_p}_{\ST_L})$;

\item {\bbf the two-level correspondences $\LGGC_{\ST-\MG}$ are related to singular universal mathematical and physical $\GL_2$-structures\/} given by the bisemisheaves 
$(\wt M^{T_p-S_p}_{\ST_R}\oplus\wt M^{T_p-S_p}_{\MG_R})\otimes
(\wt M^{T_p-S_p}_{\ST_L}\oplus\wt M^{T_p-S_p}_{\MG_L})$;

\item {\bbf the three-level correspondences $\LGGC_{\ST-\MG-\M}$ are related to singular universal mathematical and physical $\GL_2$-structures\/} given by the bisemisheaves
$(\wt M^{T_p-S_p}_{\ST_R}\oplus\wt M^{T_p-S_p}_{\MG_R}\oplus\wt M^{T_p-S_p}_{\M_R})\otimes
(\wt M^{T_p-S_p}_{\ST_L}\oplus\wt M^{T_p-S_p}_{\MG_L}\oplus\wt M^{T_p-S_p}_{\M_L})$.
\Ee

\section{Nonsingular universal $GL(2)$-structures}

\subsection{Archimedean places}

Let $k$ be a global number field of characteristic $0$.

{\bbf Let $\wt L=\wt L_R\cup\wt L_L$- denote a finite symmetric complex splitting field of $k$\/} composed of right and left algebraic extension semifields $\wt L_R$ and $\wt L_L$ in one-to-one correspondence in such a way that $\wt L_L$ \resp{$\wt L_R$} is the set of complex \resp{conjugate complex} simple roots.

Similarly, {\bbf let $\wt L^+=\wt L_R^+\cup\wt L_L^+$ be a finite symmetric real splitting field of $k$\/} where $\wt L_L^+$ \resp{$\wt L_R^+$} is the \lr algebraic extension semifield composed of the set of positive \resp{symmetric  negative} simple real roots.

Let $\wt L_{\omega _1}\subset \dots \subset \wt L_{\omega _\mu }\subset\dots \subset \wt L_{\omega _t }$
\resp{$\wt L_{\o\omega _1}\subset \dots \subset \wt L_{\o\omega _\mu }\subset\dots \subset \wt L_{\o\omega _t }$}
denote {\bbf the set of increasing complex subsemifields of $\wt L_L$
\resp{$\wt L_R$}} and let
$\wt L_{v _1}\subset \dots \subset \wt L_{v _\mu }\subset\dots \subset \wt L_{v _t }$
\resp{$\wt L_{\o v _1}\subset \dots \subset \wt L_{\o v _\mu }\subset\dots \subset \wt L_{\o v _t }$} be the {\bbf set of increasing real subsemifields of $\wt L_L^+$
\resp{$\wt L_R^+$}}.

To each complex extension 
$\wt L_{\omega _\mu }$
\resp{$\wt L_{\o\omega _\mu }$}
is associated the complex completion $ L_{\omega _\mu }$
\resp{$ L_{\o\omega _\mu }$} obtained by {\bbf an isomorphism of compactification}
of $\wt L_{\omega _\mu }$
\resp{$\wt L_{\o\omega _\mu }$} onto a closed compact subset of the complex numbers $\cit$ \resp{their conjugate complex $\cit^*$}.

Similarly, a real completion
$L_{v _\mu }$
\resp{$L_{\o v _\mu }$}
can be  obtained from the real extension $\wt L_{v _\mu }$
\resp{$\wt L_{\o v _\mu }$}
by {\bbf an isomorphism of compactification}
of $\wt L_{v _\mu }$
\resp{$\wt L_{\o v _\mu }$} onto a closed compact subset of $\rit^+$ \resp{$\rit^-$}.

Let $L_\omega =\{L_{\omega _\mu }\}_\mu $
\resp{$L_{\o\omega }=\{L_{\o\omega _\mu }\}_\mu $} denote the set of \lr complex completions (with multiplicity equal to $1$) associated with
$\wt L_L$
\resp{$\wt L_R$} and let 
$L_v=\{L_{v_{\mu ,m_\mu }}\}_{\mu ,m_\mu }$
\resp{$L_{\o v}=\{L_{\o v_{\mu ,m_\mu }}\}_{\mu ,m_\mu }$}
be the set of \lr real equivalent completions (with multiplicity $m_\mu $) associated with $\wt L_L^+$
\resp{$\wt L_R^+$}.

Each real \lr pseudoramified completion (as well as the equivalent extension) is characterized by a degree
\[
[L_{v_{\mu ,m_\mu }}:k]=*+\mu \cdot N
\rresp{[L_{\o v_{\mu ,m_\mu }}:k]=*+\mu \cdot N}
\]
which is an integer modulo $N$ where $*$ denotes an integer inferior to $N$, and where $N$ is responsible for the pseudoramification \cite{Pie1}.

Each complex \lr pseudoramified completion (as well as the equivalent extension) is characterized by a degree:
\[
[L_{\omega _{\mu }}:k]=*+\mu \cdot N\cdot m^{(\mu )}
\rresp{[L_{\o\omega _{\mu  }}:k]=*+\mu \cdot N\cdot m^{(\mu )}}
\]
where $m^{(\mu )}$ denotes the multiplicity of the $\mu $-th equivalent real completions $L_{v_{\mu ,m_\mu }}$ covering the complex corresponding completion $L_{\omega _\mu }$ according to \cite{Pie1}.

The set
$\{L_{v_{\mu ,m_\mu }}\}_{m_\mu }$
\resp{$\{L_{\o v_{\mu ,m_\mu }}\}_{m_\mu }$}
of real \lr pseudoramified equivalent completions define the {\bbf infinite pseudoramified archimedean \lr real place $v_\mu $ \resp{$\o v_\mu }$} and the set
$\{L_{\omega _{\mu ,m_\mu }}\}_{m_\mu }$
\resp{$\{L_{\o \omega _{\mu ,m_\mu }}\}_{m_\mu }$}
of complex \lr pseudoramified equivalent completions define {\bbf the infinite pseudoramified archimedean \lr complex space $\omega _\mu $ \resp{$\o\omega _\mu $}.}
\vskip 11pt

\subsection{Algebraic bilinear semigroups}

Let 
\[
W_{\wt L_\omega }=\{\Gal(\dot{\wt L_{\omega _\mu }}/k)\}^t_{\mu =1}\qquad
\rresp{W_{\wt L_{\o \omega} }=\{\Gal(\dot{\wt L_{\o \omega _\mu }}/k)\}^t_{\mu =1}}
\]
be {\bbf the global Weil group} defined as the Galois subgroup of the complex pseudoramified extensions 
$\dot{\wt L}_{\omega _\mu }$
\resp{$\dot{\wt L}_{\o\omega _\mu }$}
 characterized by extension degrees $d=0\mod N$
 and let
 $W_{\wt L_v }=\{\Gal(\dot{\wt L}_{v _\mu }/k)\}^t_{\mu =1}$
\resp{$W_{\wt L_{\o v }}=\{\Gal(\dot{\wt L}_{\o v _\mu }/k)\}^t_{\mu =1}$}
be the corresponding global Weil group of the real pseudoramified extensions
$\dot{\wt L}_{v _\mu } $
\resp{$\dot{\wt L}_{\o v _\mu } $}.

In the context of bisemistructures \cite{Pie4}, let
$G^{(2)}(L_{\o v}\times L_v)$ be the $2$-dimensional real compactified representation space of
$(W_{\wt L_{\o v}}\times W_{\wt L_v})$ covering in the sense of \cite{Pie5} its complex equivalent $G^{(1)}(L_{\o \omega }\times L_\omega )$
being the complex compactified representation space of
$(W_{\wt L_{\o \omega }}\times W_{\wt L_\omega })$,
$G^{(2)}(L_{\o v}\times L_v)$ being identified with
$G^{(1)}(L_{\o \omega }\times L_\omega )$ \cite{Del} or with
$G^{(2)}(L_{\o \omega }\times L_\omega )$ (particular ``restricted'' case considered in \cite{Pie5}).

So, {\bbf the algebraic representation of the algebraic bilinear semigroup of matrices}
$\GL_2(L_{\o v}\times L_v)$ in the
$\GL_2(L_{\o v}\times L_v)$-bisemimodule $G^{(2)}(L_{\o v}\times L_v)$, also noted
$M_R(L_{\o v})\otimes M_L(L_v)$, results from the morphism of
$\GL_2(L_{\o v}\times L_v)$ into the group of automorphisms
$\GL(M_R(L_{\o v})\otimes M_L(L_v))$ of
$M_R(L_{\o v})\otimes M_L(L_v)$.

{\bbf The algebraic bilinear semigroup $\GL_2(L_{\o v}\times L_v)$} over products of symmetric pairs of real completions
{\bbf covers its classical linear equivalent} as proved in \cite{Pie5} and has {\bbf the structure of a bisemigroup given by the triple
$(G_L,G_R,G\RL)$} \cite{Pie4} where:
\Bean
\item $G_L$ \resp{$G_R$} is a \lr semigroup under the addition of its \lr elements 
$g_{L_i}$ (\resp{$g_{R_i}$} restricted to (or referring to) the upper \resp{lower} half space.

\item $G\RL$ is a bilinear semigroup whose bielements $(g_{R_i}\times g_{L_i})$ are submitted to the cross binary operation $\u\times$ according to:
\begin{eqnarray*}
G\RL \Times G\RL &\To& G\RL\\
(g_{R_i}\times g_{L_i})\Times(g_{R_j}\times g_{L_j})&\To & 
(g_{R_i}+g_{R_j})\times (g_{L_i}+g_{L_j})
\end{eqnarray*}
leading to the cross products 
$(g_{R_i}\times g_{l_j})$ and $ (g_{R_j}\times g_{L_i})$.
\Ee

The algebraic bilinear semigroup
$\GL_2(L_{\o v}\times L_v)$ decomposes according to:
\[
\GL_2(L_{\o v}\times L_v)=T_2^t(L_{\o v})\times T_2(L_v)\]
where $T_2^t(L_{\o v})$
\resp{$T_2(L_{v})$} is a \rl linear semigroup of lower \resp{upper} triangular matrices which entries in $L_{\o v}$ \resp{$L_v$}.

The algebraic representation space $M_R(L_{\o v})\times M_L(L_v)$ of $\GL_2(L_{\o v}\times L_v)$ decomposes into a set of conjugacy class representatives on the real pseudoramified completions 
$L_{\o v_{\mu ,m_\mu }}$ and
$L_{v_{\mu ,m_\mu }}$ according to:
\[ M_R(L_{\o v})\otimes M_L(L_v)=\{M_{\o v_{\mu ,m_\mu }}\otimes
M_{v_{\mu ,m_\mu }}\}_{v_{\mu ,m_\mu }}\]
where $m_\mu $ labels the multiplicity of the $\mu $-th conjugacy class representative
$M_{\o v_{\mu ,m_\mu }}\otimes M_{v_{\mu ,m_\mu }}$.

Let $\phi _R(M_{\o v_{\mu ,m_\mu }})$
\resp{$\phi _L(M_{v_{\mu ,m_\mu }})$} be a complex-valued one-dimensional differentiable function over the
$(\mu ,m_\mu )$-th conjugacy class representative of
$T^t_2(L_{\o v})$
\resp{$T_2(L_{v})$} and let
$\phi _R(M_{\o v_{\mu ,m_\mu }})\otimes \phi _L(M_{v_{\mu ,m_\mu }})$ denote the corresponding bifunction on $M_{\o v_{\mu ,m_\mu }}\otimes M_{v_{\mu ,m_\mu }}$.

Then, the set $\{\phi _R(M_{\o v_{\mu ,m_\mu }})\}_{\mu ,m_\mu }$
\resp{$\{\phi _L(M_{v_{\mu ,m_\mu }})\}_{\mu ,m_\mu }$}
of $\cit$-valued differentiable functions, localized in the lower \resp{upper} half space is the set
$\Gamma (\phi _R(M_R(L_{\o v})))$
\resp{$\Gamma (\phi _L(M_L(L_{v})))$}
 of right \resp{left} sections of the semisheaf of rings 
 $\phi _R(M_R(L_{\o v}))$
\resp{$\phi _L(M_L(L_{v}))$} \cite{Mum}, \cite{Ser}.

And, {\bbf the set 
$\{\phi _R(M_{\o v_{\mu ,m_\mu }})\otimes \phi _L(M_{v_{\mu ,m_\mu }})\}_{\mu ,m_\mu }$
of differentiable bifunctions constitutes the set
$\Gamma (\phi _R(M_R(L_{\o v}))\otimes \phi _L(M_L(L_{v})))$ of one-dimensional bisections of the bisemisheaf of rings
$\phi _R(M_R(L_{\o v})\otimes \phi _L(M_L(L_{v})$
over $\GL_2(L_{\o v}\times L_v)$\/}.

Similarly, the set $\{\phi _R(M_{\o\omega _\mu })\otimes\phi _L(M_{\omega _\mu })\}_\mu $ of two-dimensional differentiable bifunctions constitutes the set
$\Gamma (\phi _R(M_R(L_{\o\omega}))\otimes \phi _L(M_L(L_{\omega})))$ of complex bisections of the bisemisheaf of rings
$\phi _R(M_R(L_{\o\omega}))\otimes \phi _L(M_L(L_{\omega}))$ over $\GL_1(L_{\o\omega }\times L_\omega )$.
\vskip 11pt

\subsection{Proposition (Langlands global correspondence: two-dimensional real case)}

{\em
Let
\begin{align*}
\sigma _{\mu ,m_\mu }(W_{\wt L_{\o v_{\mu ,m_\mu }}}\times W_{\wt L_{v_{\mu ,m_\mu }}})
&= G^{(2)}(L_{\o v_{\mu ,m_\mu }}\times L_{v_{\mu ,m_\mu }})\\
&= M_{\o v_{\mu ,m_\mu }}\otimes M_{\o v_{\mu ,m_\mu }}
\end{align*}
be the {\bbf $2$-dimensional representation subspace 
$(M_{\o v_{\mu ,m_\mu }}\otimes M_{\o v_{\mu ,m_\mu }})$ of the product, right by left, of the global Weil subgroups\/}
$W_{\wt L_{\o v_{\mu ,m_\mu }}}\times W_{\wt L_{v_{\mu ,m_\mu }}}$
restricted to the $\mu $-th real extensions.

Let
\[
\Pi _{\mu ,m_\mu }(\GL_2 ( L_{\o v_{\mu ,m_\mu }} \times L_{v_{\mu ,m_\mu }}))
=\Pi ^\vee _{\mu ,m_\mu }(T^t_2(L_{\o v_{\mu ,m_\mu }}))\times
 \Pi _{\mu ,m_\mu }(T_2(L_{v_{\mu ,m_\mu }}))\]
 be the {\bbf cuspidal representation of the algebraic bilinear subsemigroup}
 $\GL_2(L_{\o v_{\mu ,m_\mu }}\times L_{v_{\mu ,m_\mu }})$ in such a way that
 $\Pi ^\vee_{\mu ,m_\mu }(T^t_2(L_{\o v_{\mu ,m_\mu }}))$ be the contragradient cuspidal subrepresentation of
 $T^t_2(L_{\o v_{\mu ,m_\mu }})$ with respect to
 $T_2(L_{v_{\mu ,m_\mu }})$.
 
 Then, there exists the bijective morphism:
 \[
 \LGC_{\mu ,m_\mu }: \quad
 \sigma _{\mu ,m_\mu }(W_{\wt L_{\o v_{\mu ,m_\mu }}}\times W_{\wt L_{v_{\mu ,m_\mu }}})
 \To
 \Pi _{\mu ,m_\mu }(\GL_2 ( L_{\o v_{\mu ,m_\mu }} \times L_{v_{\mu ,m_\mu }}))
 \]
 between the $(\mu ,m_\mu )$-th representative subspace of the product, right by left, of Weil global subgroups given by the
$\GL_2 ( L_{\o v_{\mu ,m_\mu }} \times L_{v_{\mu ,m_\mu }})$-subbisemimodule
$(M_{\o v_{\mu ,m_\mu }}\otimes M_{\o v_{\mu ,m_\mu }})$ and the corresponding cuspidal class representative
$\Pi _{\mu ,m_\mu }(\GL_2 ( L_{\o v_{\mu ,m_\mu }} \times L_{v_{\mu ,m_\mu }}))
$.

This leads to the {\bbf Langlands global correspondence\/}:
\[ \LGC:\quad
\sigma (W_{\wt L_{\o v}}\times W_{\wt L_{v}})
 \To
 \Pi (\GL_2 ( L_{\o v} \times L_{v}))
 \]
 between the set
 $\sigma (W_{\wt L_{\o v}}\times W_{\wt L_{v}})$ of $2$-dimensional representative subspaces of the Weil global subgroups given by the algebraic bilinear semigroup
 $(M_R(L_{\o v})\otimes M_L(L_v))=G^{(2)}(L_{\o v}\times L_v)$ and its cuspidal representation given by $\Pi (\GL_2 ( L_{\o v} \times L_{v}))$ \cite{Pie6}.
 }
 \vskip 11pt
 
\begin{proof}
Let
\[
\TC_{L_{\mu ,m_\mu }}: \quad \wt L_{v_{\mu ,m_\mu }}\To L^T_{v_{\mu ,m_\mu }}\qquad
\rresp{\TC_{R_{\mu ,m_\mu }}: \quad \wt L_{\o v_{\mu ,m_\mu }}\To L^T_{\o v_{\mu ,m_\mu }}}\]
be the {\bbf toroidal compactification of the extension
$\wt L_{v_{\mu ,m_\mu }}$
\resp{$\wt L_{\o v_{\mu ,m_\mu }}$}
into a semicircle
$\wt L^T_{v_{\mu ,m_\mu }}$
\resp{$\wt L^T_{\o v_{\mu ,m_\mu }}$}}
localized in the upper \resp{lower} half space.

Then, by this mapping
$\TC_{L_{\mu ,m_\mu }}$
\resp{$\TC_{R_{\mu ,m_\mu }}$}, the complex-valued differentiable function
$\phi _L(M_{v_{\mu,m_\mu  }})$
\resp{$\phi _R(M_{\o v_{\mu ,m_\mu  }})$}
over $T_2(L_{v_{\mu ,m_\mu }})$
\resp{$T^t_2(L_{\o v_{\mu ,m_\mu }})$} is transformed into:
\begin{align*}
\TC_{L_{\mu ,m_\mu }}: \qquad \phi _L(M_{v_{\mu ,m_\mu }})&\To\phi _L(M_{v^T_{\mu ,m_\mu }})\\
\rresp{\TC_{R_{\mu ,m_\mu }}: \qquad \phi _R(M_{v_{\mu ,m_\mu }})&\To\phi _R(M_{\o v^T_{\mu ,m_\mu }})}
\end{align*}
in such a way that:
\[
\Pi _{\mu ,m_\mu }(T_2(L_{v_{\mu ,m_\mu }}))=\phi _L(M_{v^T_{\mu ,m_\mu }})\quad
\rresp{\Pi^\vee _{\mu ,m_\mu }(T^t_2(L_{\o v_{\mu ,m_\mu }}))=\phi _R(M_{\o v^T_{\mu ,m_\mu }})}
\]
implying that the differentiable function
$\phi _L(M_{v^T_{\mu ,m_\mu }})$
\resp{$\phi _R(M_{\o v^T_{\mu ,m_\mu }})$}
over the $(\mu ,m_\mu )$-th toroidal compactified conjugacy class representative
$T_2(L^T_{v_{\mu ,m_\mu }})$
\resp{$T^t_2(L^T_{v_{\mu ,m_\mu }})$}
of the algebraic linear semigroup
$T_2(L^T_{v})$
\resp{$T^t_2(L^T_{\o v})$}
is the cuspidal representation $\Pi _{\mu ,m_\mu }(T_2(L_{v_{\mu ,m_\mu }}))$
\resp{contragradient cuspidal representation $\Pi^\vee _{\mu ,m_\mu }(T^t_2(L_{\o v_{\mu ,m_\mu }}))$}
of $T_2(L_{v_{\mu ,m_\mu }})$
\resp{\linebreak $T^t_2(L_{\o v_{\mu ,m_\mu }})$},
where $L^T_v=\{L^T_{v_{\mu ,m_\mu }}\}_{\mu ,m_\mu }$.

Indeed, by summing over all $\mu $ and $m_\mu $, we get the global elliptic semimodule \cite{Pie5}:
\begin{align*}
\phi _L(M_L(L^T_v)) &= \sum_{\mu ,m_\mu }\phi _L(M_{v^T_{\mu ,m_\mu }})\\
 &= \sum_{\mu ,m_\mu }r(\mu ,m_\mu )\ e^{\pi i\mu x_L}\;, && x_L\in L_v\;, \\
\rresp{\phi _R(M_R(L^T_{\o v})) &= \sum_{\mu ,m_\mu }\phi _R(M_{\o v^T_{\mu ,m_\mu }})\\
 &= \sum_{\mu ,m_\mu }r(\mu ,m_\mu )\ e^{-\pi i\mu x_L}}\;,
 \end{align*}
 where 
 $r(\mu ,m_\mu )\ e^{\pi i\mu x_L}$
\resp{ $r(\mu ,m_\mu )\ e^{-\pi i\mu x_L}$}
is a semicircle localized in the upper \resp{lower} half plane.

According to \cite{Pie5}, the {\bbf global elliptic semimodule
$\phi _L(M_L(L^T_v))$
\resp{$\phi _R(M_R(L^T_{\o v}))$}
over 
$T_2(L^T_{v})$
\resp{$T^t_2(L^T_{\o v})$}
covers (and is thus isomorphic to) the corresponding cuspidal form over
$T_1(L_\omega )$
\resp{$T^t_1(L_{\o \omega} )$}}, where
$GL_1(L_{\o\omega} \times L_\omega )=T^t_1(L_{\o\omega })\times T_1(L_\omega )$.

Thus, $\phi _L(M_{v^T_{\mu ,m_\mu }})$
\resp{$\phi _R(M_{\o v^T_{\mu ,m_\mu }})$}
is a cuspidal representation restricted to the
$(\mu ,m_\mu )$-th conjugacy class representative of
$T_2(L_{v})$ \resp{$T^t_2(L_{\o v })$}.

So, by considering the toroidal compactification of all $(\mu ,m_\mu )$ extensions or completions, we get the Langlands bilinear global correspondence
\[\LGC: \qquad \sigma (W_{\wt L_{\o v}}\times W_{\wt L_v})\To \Pi (\GL_2(L_{\o v}\times L_v))\]
equivalent to the isomorphism:
\[
\LGC: \qquad
\phi _R(M_R(L_{\o v}))\otimes\phi _L(M_L(L_{v}))
\To \Pi (\GL_2(L_{\o v}\times L_v))\]
where $\phi _R(M_R(L_{\o v}))\otimes\phi _L(M_L(L_{v}))$, which is a bisemisheaf of rings over the algebraic bilinear semigroup
$\GL_2(L_{\o v}\times L_v)$, constitutes the two dimensional functional representation space 
$\sigma (W_{\wt L_{\o v}}\times W_{\wt L_v})$ of the product, right by left, of global Weil groups.
\end{proof}
\vskip 11pt

\subsection{Proposition}

{\em The bisemisheaf of rings $\phi _R(M_R(L_{\o v})\otimes \phi _L(M_L(L_{v})$ is a physical bosonic quantum string field of an elementary particle\/}.
\vskip 11pt

\begin{proof} Let $L^T_{\o v}=\{L^T_{\o v_{\mu ,m_\mu }}\}_{\mu ,m_\mu }$
\resp{$L^T_{v}=\{L^T_{v_{\mu ,m_\mu }}\}_{\mu ,m_\mu }$} be the set  of \rl real pseudoramified {\bf toroidal completions\/} obtained from $L_{\o v}$ \resp{$L_v$} by a toroidal compactification of the corresponding completions
$L_{\o v_{\mu ,m_\mu }}$
\resp{$L_{v_{\mu ,m_\mu }}$} \cite{Pie3}.

The set $\{L^T_{\o v_{\mu ,m_\mu }}\}_{m_\mu }$
\resp{$\{L^T_{v_{\mu ,m_\mu }}\}_{m_\mu }$}
of $\mu $-th completions are then semicircles covering a $2$-dimensional \rl semitorus $T^2_R(\mu )$ \resp{$T^2_L(\mu )$}, localized in the lower \resp{upper} half space.
\vskip 11pt

Assume that the degree of $L^T_{\o v_{\mu ,m_\mu }}$ and of
$L^T_{v_{\mu ,m_\mu }}$ is given by
\[[L^T_{\o v_{\mu ,m_\mu }}]\equiv [L^T_{v_{\mu ,m_\mu }}]=\mu \ N\;.\]
Then, $L^T_{\o v_{\mu ,m_\mu }}$
and $L^T_{v_{\mu ,m_\mu }}$ are toroidal completions or semicircles at $\mu $ quanta, a quantum being an irreducible completion of degree $N$.
\vskip 11pt

As we are essentially interested in circles, toroidal completions
$L^T_{\o v_{2\mu ,m_{2\mu }}}$
and $L^T_{v_{2\mu ,m_{2\mu} }}$ characterized by degrees:
\[[L^T_{\o v_{2\mu ,m_{2\mu} }}]\equiv [L^T_{v_{2\mu ,m_{2\mu} }}]=2\ \mu \ N\]
will be taken into account.

By this way, the corresponding completions $L_{\o v_{2\mu ,m_{2\mu }}}$
and $L_{v_{2\mu ,m_{2\mu} }}$ are closed paths or closed strings.
\vskip 11pt

{Now, each product
$\{L^T_{\o v_{2\mu ,m_{2\mu }}}\otimes L^T_{v_{2\mu ,m_{2\mu} }}\}$ of symmetric circles\/} rotating in opposite senses according to \cite{Pie2} {is the representation of an harmonic oscillator\/}.  This is also the case for the product $\{L_{\o v_{2\mu ,m_{2\mu }}}\otimes L_{v_{2\mu ,m_{2\mu} }}\}$ of corresponding completions homeomorphic to
$\{L^T_{\o v_{2\mu ,m_{2\mu }}}\otimes L^T_{v_{2\mu ,m_{2\mu} }}\}$ and for the product
$\{\phi _R(M_{\o v_{2\mu ,m_{2\mu }}})\otimes \phi _L(M_{v_{2\mu ,m_{2\mu} }})\}$ of $\cit$-valued differentiable functions.
\vskip 11pt

Consequently, {the set of packets
$\{\phi _R(M_{\o v_{2\mu ,m_{2\mu }}})\otimes \phi _L(M_{v_{2\mu ,m_{2\mu} }})\}_{2\mu ,m_{2\mu }}$ of bifunctions\/} on\linebreak
$\{L{\o v_{2\mu ,m_{2\mu }}}\otimes L_{v_{2\mu ,m_{2\mu} }}\}_{2\mu ,m_{2\mu }}$ {behaves like a set of packets of harmonic oscillators characterized by increasing integers $2\mu $\/}.

Thus, the bisemisheaf of rings
\[\phi _R(M_R(L_{\o v}))\otimes \phi _L(M_L(L_{v}))=
\{\phi _R(M_{\o v_{2\mu ,m_{2\mu }}})\otimes \phi _L(M_{v_{2\mu ,m_{2\mu} }})\}_{2\mu ,m_{2\mu }}\] is a physical string field.
\vskip 11pt

{It is a quantum string field because the set of sections of the bisemisheaf
$\phi _R(M_R(L_{\o v}))\linebreak\otimes \phi _L(M_L(L_{v}))$ is a tower of increasing bistrings\/}, i.e. products of symmetric right and left strings,  behaving like harmonic oscillators and characterized by a number of increasing biquanta, $2\mu \le\infty $, corresponding to the normal modes of the string field.

{This quantum string field is a bosonic field\/} because biquanta can be added to (i.e. created) or removed (i.e. annihilated) from these bistrings by Galois automorphisms or antiautomorphisms as developed in \cite{Pie1} and because each bistring with degree $2\mu $ was interpreted as a ``bound bisemiphoton'' at $2\mu $ biquanta.
\end{proof}
\vskip 11pt

\subsection{Properties of the bisemisheaf $\phi _R(M_R(L_{\o v}))\otimes \phi _L(M_L(L_{v}))$}

The representation semispace 
$\Repsp (T_2(L_v))$
\resp{$\Repsp (T^t_2(L_{\o v}))$} of the linear algebraic semigroup
$T_2(L_{v})$
\resp{$T^t_2(L_{\o v})$} is the unitary \lr 
$L_v$-semimodule $M_L(L_v)$
\resp{$L_{\o v}$-semimodule $M_R(L_{\o v})$}, i.e. 
a left vector $L_v$-semispace
\resp{a right vector $L_{\o v}$-semispace}.

And the \lr semisheaf $\phi _L(M_L(L_{v}))$
\resp{$\phi _R(M_R(L_{\o v}))$}
is also a \lr vector $L_v$-semispace
\resp{$L_{\o v}$-semispace}
implying that the bisemisheaf
$\phi _R(M_R(L_{\o v}))\otimes \phi _L(M_L(L_{v}))$ is a vector
$L_{\o v}\times L_v$-bisemispace as developed in \cite{Pie4}.
\vskip 11pt

{\bbf This vector $L_{\o v}\times L_v$-bisemispace
$\phi _R(M_R(L_{\o v}))\otimes_{L_{\o v}\times L_v}\phi _L(M_L(L_{v}))$
splits naturally into\/}:
\begin{multline*}
\phi _R(M_R(L_{\o v}))\otimes_{L_{\o v}\times L_v}\phi _L(M_L(L_{v}))\\
=(\phi _R(M_R(L_{\o v}))\otimes_{D}\phi _L(M_L(L_{v}))\oplus
(\phi _R(M_R(L_{\o v}))\otimes_{OD}\phi _L(M_L(L_{v})))
\end{multline*}
where:
\Bi
\item $\phi _R(M_R(L_{\o v}))\otimes_{D}\phi _L(M_L(L_{v}))$ is a diagonal vector
$L_{\o v}\times L_{v}$-bisemispace characterized by a bilinear diagonal basis
$\{e_\alpha \otimes f_\alpha \}^2_{\alpha =1}$ of dimension $2$;

\item $(\phi _R(M_R(L_{\o v}))\otimes_{OD}\phi _L(M_L(L_{v}))$ is an off-diagonal vector
$L_{\o v}\times L_{v}$-bisemispace of dimension $2$ characterized by a bilinear off-diagonal basis
$\{e_\alpha \otimes f_\beta  \}^2_{\alpha \neq \beta =1}$.
\Ei
\vskip 11pt

As it was seen in \cite{Pie4}, the vector 
$L_{\o v}\times L_{v}$-bisemispace
$\phi _R(M_R(L_{\o v}))\otimes_{D}\phi _L(M_L(L_{v}))$, endowed with a suitable  inner product at the condition that the \rl vector $L_{\o v}$-semispace
$\phi _R(M_R(L_{\o v}))$ \resp{$\phi _L(M_L(L_{v}))$} be projected onto its symmetric \lr equivalent
$\phi _L(M_L(L_{v}))$
\resp{$\phi _R(M_R(L_{\o v}))$}, can give rise to an inner product bisemispace.
\vskip 11pt

Consequently, {\bbf this inner product bisemispace can generate an orthogonal complement bisemispace} as it is developed in the next section.
\vskip 11pt

\subsection{Endomorphisms $E_L$ and $E_R$ based on Galois antiautomorphisms}

The representation semispace 
$M_L(L_v)=\Repsp(T_2(L_v))$
\resp{$M_R(L_{\o v})=\Repsp(T^t_2(L_{\o v}))$}
of the linear \lr algebraic semigroup
$T_2(L_v)$
\resp{$T^t_2(L_{\o v})$} is noetherian or solvable in the sense that it is composed of the set
\begin{align*}
&M_L(L_{v_1})\subset \dots \subset M_L(L_{v_\mu })\subset \dots \subset M_L(L_{v_t })\\
\rresp{&M_R(L_{\o v_1})\subset \dots \subset M_R(L_{\o v_\mu })\subset \dots \subset M_R(L_{\o v_t })}
\end{align*}
of embedded increasing representation subsemispaces.

So, {\bbf we can define the smooth endomorphism\/}:
\begin{align*}
E_L: \qquad M_L(L_v) &\overset{\sim}{\To} M_L^{(r)}(L_v)\oplus M_L^{(I)}(L_v)\\
\rresp{E_R: \qquad M_R(L_{\o v}) &\overset{\sim}{\To} M_R^{(r)}(L_{\o v})\oplus M_R^{(I)}(L_{\o v})}
\end{align*}
decomposing $M_L(L_v)$
\resp{$M_R(L_{\o v})$} into the direct sum of the  reduced representation semispace
$M_L^{(r)}(L_v)$
\resp{$M_R^{(r)}(L_{\o v})$},
{\bbf submitted to Galois antiautomorphisms\/}, and of the complementary representation semispace
$M_L^{(I)}(L_v)$
\resp{$M_R^{(I)}(L_{\o v})$}, submitted to Galois automorphisms.

Similarly, the semisheaf
$\phi _L(M_L(L_v))$
\resp{$\phi _R(M_L(L_{\o v}))$}
can be submitted to the same endomorphism $E_L$ \resp{$E_R$} transforming it into:
\begin{align*}
E_L(\phi _L(M_L(L_v))) &= \phi _L(M^{(r)}_L(L_v))\oplus\phi _L(M^{(I)}_L(L_v))\\
\rresp{E_R(\phi _R(M_R(L_{\o v}))) &= \phi _R(M^{(r)}_R(L_{\o v}))\oplus\phi _R(M^{(I)}_R(L_{\o v}))}.
\end{align*}
\vskip 11pt

\subsection{Proposition (Generation of an orthogonal complement (semi)sheaf)}

{\em The semisheaf
$\phi _L(M_L(L_v))$
\resp{$\phi _R(M_R(L_{\o v}))$}
on the representation semispace
$M_L(L_v)=\Repsp(T_2(L_v))$
\resp{$M_R(L_{\o v})=\Repsp(T^t_2(L_{\o v}))$}
of the algebraic semigroup
$T_2(L_v)$
\resp{$T^t_2(L_{\o v})$}
can generate {\bbf the orthogonal complement semisheaf}
$\phi _L^\perp(\Repsp T_2(L_v))$
\resp{\linebreak $\phi _R^\perp(\Repsp T^t_2(L_{\o v}))$}
{\bbf under the composition of morphisms\/}
\begin{align*}
\gamma _{t\to r}\circ E_L:\qquad
\phi _L(M_L(L_v))&\To \phi _L(M^{(r)}_L(L_v))\oplus\phi ^\perp_L(\Repsp T_2(L_v))\\
\rresp{\gamma _{t\to r}\circ E_R:\qquad
\phi _R(M_R(L_{\o v}))&\To \phi _R(M^{(r)}_R(L_{\o v}))\oplus\phi ^\perp_R(\Repsp T^t_2(L_{\o v})}
\end{align*}
where $\gamma _{t\to r}$ {\bbf is the emergent morphism \cite{Pie8} mapping the complementary semisheaf
$\phi^{(I)} _L(M_L(L_v))$
\resp{$\phi^{(I)} _R(M_R(L_{\o v}))$}
throughout the origin into the orthogonal complement semisheaf
$\phi ^\perp_L(\Repsp T_2(L_v))$
\resp{$\phi ^\perp_R(\Repsp T^t_2(L_{\o v}))$}.
}}
\vskip 11pt

\begin{proof}
In fact, the endomorphism $E_L$ \resp{$E_R$}, acting by means of Galois antiautomorphisms, generates the reduced semisheaf
$\phi _L(M_L^{(r)}(L_v))$
\resp{$\phi _R(M_R^{(r)}(L_{\o v}))$}
and the complementary disconnected semisheaf
$\phi _L(M_L^{(I)}(L_v))$
\resp{$\phi _R(M_R^{(I)}(L_{\o v}))$}
as developed in \cite{Pie8}.

As 
$(\phi _R(M_R(L_{\o v}))\otimes\phi _L(M_L(L_{v})))$ can give rise to an inner product bisemispace \cite{Pie4} due to\linebreak the symmetry between the right and left semisheaves
$\phi _R(M_R(L_{\o v}))$ and $\phi _L(M_L(L_{v}))$,
the disconnected complementary semisheaf
$\phi ^{(I)}_L(M_L(L_{v}))$
\resp{$\phi ^{(I)}_R(M_R(L_{{\o v}}))$} can be projected by the emergent morphism
$\gamma _{t\to r}$ into an orthogonal complement (semi)space throughout the origin,\linebreak  generating then the orthogonal complement semisheaf
$\phi _L^\perp(\Repsp(T_2(L_v))$
\resp{\linebreak $\phi _R^\perp(\Repsp(T^t_2(L_{\o v}))$}.
\end{proof}
\vskip 11pt

\subsection{Proposition (Langlands general global correspondence)}

{\em 
Taking into account the existence of an orthogonal complement bisemispace, the following Langlands general global correspondence can be stated:
\[ \LGGC_\ST: \qquad
\sigma (W_{\wt L_{\o v}}\times W_{\wt L_v})
\overset{\sim}{\To}
\Pi (\GL_2^{(r)}(L_{\o v}\times L_v))+
\Pi (\GL_2^\perp(L_{\o v}\times L_v))\]
where:
\Bi
\item $\Pi (\GL_2^{(r)}(L_{\o v}\times L_v))$ is the cuspidal representation on the reduced algebraic bilinear semigroup
$\GL_2^{(r)}(L_{\o v}\times L_v)$;

\item $\Pi (\GL_2^\perp(L_{\o v}\times L_v))$ is the cuspidal representation on the orthogonal complement algebraic bilinear semigroup
$\GL_2^\perp(L_{\o v}\times L_v)$.
\Ei}
\vskip 11pt

\begin{proof}
{\bbf Under the toroidal compactification $\TC$ (see proposition 2.3), the semisheaf
$\phi _L(M_L(L_v))$
\resp{$\phi _R(M_R(L_{\o v}))$}
is transformed into the cuspidal representation
$\Pi (T_2(L_v))$
\resp{$\Pi (T^t_2(L_{\o v}))$}
of the algebraic semigroup
$T_2(L_v)$
\resp{$T^t_2(L_{\o v})$}}
and the bisemisheaf
$\phi _R(M_R(L_{\o v}))\otimes\phi _L(M_L(L_{v}))$ is transformed by:
\[\TC_R\times \TC_L: \qquad
\phi _R(M_R(L_{\o v}))\otimes\phi _L(M_L(L_{v}))\To \Pi (\GL_2(L_{\o v}\times L_v))\]
into the cuspidal representation
$\Pi (\GL_2(L_{\o v}\times L_v)$ of the algebraic bilinear semigroup
$\GL_2(L_{\o v}\times L_v)$.  So, we get the commutative diagram:
\[
\begin{CD}
\LGC :\\[8pt] @VVV\\[8pt] \LGGC_\ST :\\[24pt]
\end{CD}
 \qquad 
\begin{CD}
\sigma (W_{\wt L_{\o v}})\times \sigma (W_{\wt L_{v}})
@>{\sim}>{\TC_R\times \TC_L}> \Pi (\GL_2(L_{\o v}\times L_v))\\
 @VV{\begin{smallmatrix}
(\gamma _{t\to r}\times \gamma _{t\to r})\\
\circ (E_R\times E_L) \end{smallmatrix}}V
@VV{\begin{smallmatrix}
(\gamma _{t\to r}\times \gamma _{t\to r}\\
\circ (E_R\times E_L) \end{smallmatrix}}V\\
\begin{array}[t]{l}
\sigma^{(r)} (W_{\wt L_{\o v}}\times W_{\wt L_{v}})\\[-11pt] \quad \oplus
\sigma^{\perp} (W_{\wt L_{\o v}}\times W_{\wt L_{v}})\end{array}
@>{\sim}>> 
\begin{array}[t]{l}
\Pi (\GL^{(r)}_2(L_{\o v}\times L_v))\\[-11pt] \quad \oplus
 \Pi (\GL^{\perp}_2(L_{\o v}\times L_v))\end{array}
 \end{CD}\]
 leading to the Langlands general global correspondence
 $\LGGC_{\ST}$ where
 $\sigma^{\perp} (W_{\wt L_{\o v}}\times W_{\wt L_{v}})
 = \phi _R^\perp(M_R^\perp(L_{\o v}))\oplus
  \phi _L^\perp(M_L^\perp(L_{v}))$ is the orthogonal complement representation bisemispace of the product, right by left, of global Weil groups.
 \end{proof}
 \vskip 11pt
 
 \subsection{Proposition (Nonsingular universal $\GL(2)$-structures)}
 
 {\em \Bena
 \item The functional representation bisemispace
 $\sigma (W_{\wt L_{\o v}}\times W_{\wt L_v})$ of the product, right by left, of global Weil groups, given by
 \[
\sigma (W_{\wt L_{\o v}}\times W_{\wt L_v})=
(\phi _R(M_R^{(r)}(L_{\o v}))\otimes(\phi _L(M_L^{(r)}(L_{v}))\oplus
(\phi ^\perp_R(M_R^\perp(L_{\o v}))\otimes(\phi ^\perp_L(M_L^\perp(L_{v}))
\]
the direct sum of the reduced bisemisheaf over the reduced algebraic bilinear semigroup
$\GL_2^{(r)}(L_{\o v}\times L_v)$ and of the orthogonal complement bisemisheaf over
$\GL^\perp_2(L_{\o v}\times L_v)$, is a nonsingular universal mathematical structure.

\item It is also a nonsingular universal physical structure because it corresponds to the space-time string fields of the dark energy structure of an elementary particle.
\Ee
}
\vskip 11pt

\begin{proof}
\Bena
\item The reduced bisemisheaf
$\phi _R(M_R^{(r)}(L_{\o v}))\otimes(\phi _L(M_L^{(r)}(L_{v}))$ and the orthogonal complement bisemisheaf
$\phi ^\perp_R(M_R^\perp(L_{\o v}))\otimes(\phi ^\perp_L(M_L^\perp(L_{v})$ constitute a universal mathematical structure because:
\Be
\item they are generated from the product, right by left, of global Weil groups
$(W_{\wt L_{\o v}}\times W_{\wt L_v})$;
\item they are in one-to-one correspondence with the holomorphic and cuspidal representations of
$\GL_2^{(r)}(L_{\o v}\times L_v)+\GL_2^{\perp}(L_{\o v}\times L_v)$ by means of the Langlands general global correspondence 
$\LGGC_\ST$.
\Ee

\item The reduced bisemisheaf
$\phi _R(M_R^{(r)}(L_{\o v}))\otimes(\phi _L(M_L^{(r)}(L_{v}))$, as well as the bisemisheaf\linebreak
$\phi _R(M_R(L_{\o v}))\otimes\phi _L(M_L(L_{v}))$, has one-dimensional bisections according to section 2.2.  It is thus a time string field of an elementary particle.  On the other hand, the ``diagonal'' orthogonal complement bisemisheaf
$\phi^{\perp} _R(M_R^\perp(L_{\o v}))\otimes_D(\phi^{\perp} _L(M_L^\perp(L_{v}))$, generated from the reduced bisemisheaf by the
$(\gamma _{t\to r}\times \gamma _{t\to r})\circ (E_R\times E_L)$ morphisms, may be two- or three-dimensional according to \cite{Pie8}.

Consequently, this orthogonal complement bisemisheaf corresponds to the space string field of an elementary particle generated from the reduced time string field.
\Ee

Let us abbreviate our notations:
\Bi
\item the time reduced bisemisheaf will be rewritten
\[ \wt M_{\ST_R}^T\otimes \wt M_{\ST_L}^T\;, \qquad \text{``$T$'' for time, \quad ``$\ST$'' for space-time;}\]
\item the space orthogonal complement bisemisheaf will be rewritten
\[ \wt M_{\ST_R}^S\otimes \wt M_{\ST_L}^S\;, \qquad \text{``$S$'' for space.}\]
\Ei
This space-time bisemisheaf
\[ (\wt M_{\ST_R}^T\otimes \wt M_{\ST_L}^T)\oplus
(\wt M_{\ST_R}^S\otimes \wt M_{\ST_L}^S)\]
constitutes the space-time string fields of the dark energy structure of an elementary fermion (lepton, quark or neutrino).

Indeed, the Zel'dovich's idea of connecting the vacuum energy density of quantum field theories (QFT) with the cosmological constant $\Lambda $ of general relativity (GR) led the author to propose a unification of these two theories by a new interpretation of the equations of GR being then in one-to-one correspondence with the equations of the internal dynamics of the vacuum and mass structures of a set of interacting elementary particles.  In this context, the cosmological constant would deal with the internal vacuum substructure ``$\ST$'' of the elementary particles.

The dark energy structure of an elementary particle then refers to the space-time fields of its internal vacuum structure ``$\ST$''.

As these fields are in fact bisemifields, i.e. products  of right semifields by their symmetric left correspondents, an elementary particle at this dark energy level is a bisemiparticle composed of the product of a left semiparticle, localized in the upper half space, by its symmetric right (co)semiparticle, localized in the lower half space.
\end{proof}
\vskip 11pt

\subsection{Shifted bisemisheaves}

According to proposition 2.7, the ``space'' bisemisheaf
$\wt M_{\ST_R}^S\otimes\wt M_{\ST_L}^S$ (and the ``time'' bisemisheaf
$\wt M_{\ST_R}^T\otimes\wt M_{\ST_L}^T$) constitutes the functional representation space
$\FRepsp(\GL_2(L_{\o v}\times L_v))$ of the algebraic bilinear semigroup
$\GL_2(L_{\o v}\times L_v)$.

{\bbf The dynamics of the bisemisheaf
$\wt M_{\ST_R}^S\otimes\wt M_{\ST_L}^S$ is obtained by the action of the elliptic differential bioperator
$(D_R\otimes D_L)$} mapping
$\wt M_{\ST_R}^S\otimes\wt M_{\ST_L}^S$ into its shifted equivalent according to:
\[ D_R\otimes D_L: \qquad
\wt M_{\ST_R}^S\otimes\wt M_{\ST_L}^S\To
\wt M_{\ST_R}^{S_p}\otimes\wt M_{\ST_L}^{S_p}\]
where the bisemisheaf
$\wt M_{\ST_R}^{S_p}\otimes\wt M_{\ST_L}^{S_p}$ is the functional representation space
$\FRepsp(\GL_2((L_{\o v}\otimes \rit))\times (L_v\otimes\rit))$ of the shifted bilinear algebraic semigroup 
$\GL_2((L_{\o v}\otimes \rit))\times (L_v\otimes\rit)$ as it is developed in \cite{Pie3}.

Physically, 
$(\wt M_{\ST_R}^{S_p}\otimes\wt M_{\ST_L}^{S_p})$ is an operator valued space string field if it is referred to proposition 2.9.  Similarly, the ``time'' bisemisheaf
$\wt M_{\ST_R}^{T}\otimes\wt M_{\ST_L}^{T}$ is sent by the action of an elliptic differential bioperator into its shifted equivalent
$(\wt M_{\ST_R}^{T_p}\otimes\wt M_{\ST_L}^{T_p})$ responsible for its dynamics.
\section{Singular universal $GL(2)$-structures}

In this chapter, it will be shown that {\bbf degenerate singularities on the functional representation spaces
$\FRepsp(\GL_2((L_{\o v}\otimes \rit)\times(L_v\otimes\rit))$\/}, given by the space or time bisemisheaves
$(\wt M_{\ST_R}^{S_p}\otimes \wt M_{\ST_L}^{S_p})$ or
$(\wt M_{\ST_R}^{T_p}\otimes \wt M_{\ST_L}^{T_p})$, {\bbf can give rise, by versal deformations and blowups of these, to one or two new covering functional representation spaces of $\GL_2((L_{\o v}\otimes \rit)\times(L_v\otimes\rit))$ according to the kind of considered degenerate singularities\/}.
\vskip 11pt

\subsection{One- and two-dimensional sections of the bisemisheaves}

Let $(\wt M^{S_p}_{ST_R}\otimes \wt M^{S_p}_{ST_L})$ denote the  space bisemisheaf introduced in section 2.10.
\vskip 11pt

Let $\{\phi _R^{S_p}(M_{\o v_{2\mu ,m_{2\mu }}})\otimes\phi _L^{S_p}(M_{v_{2\mu ,m_{2\mu }}})\}_{\mu ,m_\mu }$ denote the set of bisections of $(\wt M^{S_p}_{ST_R}\otimes \wt M^{S_p}_{ST_L})$
as introduced in proposition 2.4 and noted there
$\phi _R^{S_p}(M_R(L_{\o v}))\otimes \phi _L^{S_p}(M_L(L_{v}))$.

Then, 
$\{\phi _L^{S_p}(M_{v_{2\mu ,m_{2\mu }}})\}_{\mu ,m_\mu }$
\resp{$\{\phi _R^{S_p}(M_{\o v_{2\mu ,m_{2\mu }}})\}_{\mu ,m_\mu }$} will be the set
$\Gamma (\phi _L^{S_p}(M_L(L_v)))$
\resp{$\Gamma (\phi _R^{S_p}(M_R(L_{\o v})))$}
of \lr one-dimensional sections of the \lr semisheaf
$\phi _L^{S_p}(M_L(L_v))$
\resp{$\phi _R^{S_p}(M_R(L_{\o v}))$} \cite{Har}.
\vskip 11pt

Assume that the ``$m_\mu $'' one-dimensional sections
$\{\phi  _L^{S_p}(M_{v_{2\mu ,m_{2\mu }}})\}_{m_\mu }$
\resp{$\{\phi  _R^{S_p}(M_{\o v_{2\mu ,m_{2\mu }}})\}_{m_\mu }$}
of each packet ``$\mu $'' are glued together under a compactification map ``$c$'' in order to generate a surface
$\phi  _L^{S_p}(M_{v_{2\mu}(c)})$
\resp{$\phi  _R^{S_p}(M_{\o v_{2\mu}(c)})$} as described in \cite{Pie1} and in \cite{Pie7}.

Then, 
$\{\phi  _L^{S_p}(M_{v_{2\mu}(c)})\}_\mu $
\resp{$\{\phi  _R^{S_p}(M_{\o v_{2\mu}(c)})\}_\mu $} will denote the set
$\Gamma (\phi  _L^{S_p}(M_L(L_{v(c)})))$
\resp{\linebreak $\Gamma (\phi  _R^{S_p}(M_R(L_{\o v(c)})))$} of \lr two-dimensional sections of the \lr semisheaf
$\phi  _L^{S_p}(M_L(L_{v(c)}))$
\resp{$\phi  _R^{S_p}(M_R(L_{\o v(c)}))$}.
\vskip 11pt

\subsection{Degenerate singularities on the sections of the space bisemisheaves}

Under external perturbations due to the strong fluctuations at this length scale, degenerate singularities are produced on the left and right sections of the above mentioned semisheaves.  Furthermore, it is assumed that a same kind of degenerate singularities is generated on each \lr section which is a \lr differentiable function \resp{cofunction}.
\vskip 11pt

On one- or two-dimensional sections, {\bbf the simple germs $f(x)=x^{k+1}$, $1\le k\le 3$, which are singular points of corank $1$ and multiplicity ``$k-1$'', can be produced in a $3$-dimensional (semi)space by singularization morphisms\/} which are defined as sets of contracting surjective morphisms of normal crossing divisors as developed in \cite{Pie2}.
\vskip 11pt

On a two-dimensional section, {\bbf the possible degenerate singular points are essentially the following germs of corank 2 and multiplicity inferior or equal to $3$\/} \cite{Tho}, \cite{Arn}:
\[ f(x,y)=x^3-3y^2x\;, \quad f(x,y)=x^3+y^3\;.\]
They are also produced by a set of contracting surjective morphisms of normal crossing divisors.
\vskip 11pt

\subsection{Versal deformations of degenerate singularities}

Under the same kind of external perturbations, these degenerate singular germs are submitted to {\bbf versal deformations interpreted as extensions of the contracting surjective morphisms of singularizations\/} as proved in \cite{Pie2}.
\vskip 11pt

{\bbf The versal deformation or unfolding of a germ $f(x)=x^{k+1}$\/} of corank $1$ and multiplicity ``$k-1$'' is given, in the frame of the Malgrange preparation theorem, by \cite{Mal}, \cite{Mat}:
\[
F(x,a(y,z)) = f(x)+\sum_{i=1}^{k-1}a_i(y,z)\ x^i\;. \]
$R(x,a_i(y,z)) = \sum\limits_{i=1}^{k-1}a_i(y,z)\ x^i$
is the polynomial of the quotient algebra of the versal unfolding of the degenerate germ $f(x)$ with:
\Bi
\item $\{x^1,\dots,x^i,\dots,x^{k-1}\}$ being the basis of this quotient algebra of dimension $(k-1)$, which is thus finitely generated;

\item $a_i(y,z)$ being a function of two variables on which is mapped the respective generator $x^i$ as developed in \cite{Tho} and \cite{Mal}.
\Ei

This function $a_i(y,z)$ belonging to a section of the space  semisheaf is two-dimensional because the set of sections of the \lr semisheaf
$\phi _L^{S_p}(M_L(L_v))$
\resp{$\phi _R^{S_p}(M_R(L_{\o v}))$} is assumed to be compactified in a three-dimensional spatial semispace as developed in \cite{Pie8}.
\vskip 11pt

More concretely, {\bbf the versal unfolding of the degenerate germs $f(x)=x^{k+1}$ in codimension inferior or equal to thee are\/} (to a translation):
\Bena
\item {\bbf the fold\/}: \quad $f(x)=x^3$

of which versal unfolding in codimension $1$ is
\[ F(x,a_1)=x^3+a_1  x^1\;;\]

\item {\bbf the cusp\/}: \quad $f(x)=x^4$

of which versal unfolding in codimension $2$ is
\[ F(x,a_1,a_2)=x^4+a_1  x^1+a_2   x^2\;;\]

\item {\bbf the swallowtail\/}: \quad $f(x)=x^5$

of which versal unfolding in codimension $3$ is
\[ F(x,a_1,a_2,a_3)=x^5+a_1  x^1+a_2  x^2+a_3  x^3\;;\]
\Ee
where $a_i$ is a contracted notation for $a_i(y,z)$, $1\le i\le 3$.
\vskip 11pt

{\bbf The degenerate germs of corank $2$ and multiplicity inferior of equal to $3$ are\/} \cite{Tho} (to a translation):
\Bena
\item {\bbf the elliptic umbilic\/}: \quad $f(x,y)=x^3-3xy^2$

of which versal unfolding in codimension $3$ is
\[ F(x,y,b_1,b_2)=x^3-3xy^2+b_1(x^2+y^2)-b_2y\;;\]

\item {\bbf the hyperbolic umbilic\/}: \quad $f(x,y)=x^3+y^3$

of which versal unfolding in codimension $3$ is
\[ F(x,y,b_2,b_3,b_4)=x^3+y^3-b_2y-b_3x+b_4xy\]
where $b_2$ and $b_3$ are two variable functions and $b_1$ and $b_4$ are one variable functions on a section of the space  semisheaf.
\Ee
\vskip 11pt

More generally, let 
$\phi _L^{\star S_p}(M_L)$
\resp{$\phi _R^{\star S_p}(M_R)$} denote the \lr semisheaf  of which sections:
\Bena
\item are one- or two-dimensional;
\item are affected by the same kind of degenerate singularities of corank $1$ or $2$ as developed in sections 3.1 and 3.2.
\Ee

Then, {\bbf the versal deformation\/} of the semisheaf
$\phi _L^{\star S_p}(M_L)$
\resp{$\phi _R^{\star S_p}(M_R)$} of differentiable functions \resp{cofunctions} endowed with singular germs of corank $1$ or $2$ {\bbf is given by the contracting fiber bundle\/}:
\begin{align*}
D_{S_L}: \qquad  \phi _L^{\star S_p}(M_L)\times\phi _{S_L}& \To  \phi _L^{\star S_p}(M_L)\\
\rresp{D_{S_R}: \qquad  \phi _R^{\star S_p}(M_R)\times\phi _{S_R}& \To  \phi _R^{\star S_p}(M_R)}
\end{align*}
in such a way that the fiber 
$\phi _{S_L}=\{\dots,\phi _L(x^i),\dots\}_{i=1}^{k-1}$
\resp{$\phi _{S_R}=\{\dots,\phi _R(x^i),\dots\}_{i=1}^{k-1}$}
is given by the set of $(k-1)$ semisheaves of the base $S_L$ \resp{$S_R$} of the considered versal deformation where 
$\phi _L(x^i)$
\resp{$\phi _R(x^i)$} denotes the semisheaf of monomials ``$x^i$'' with respect to all the sections of
$\phi _L^{\star S_p}(M_L)$
\resp{$\phi _R^{\star S_p}(M_R)$}.

These semisheaves 
$\phi _L(x^i)$
\resp{$\phi _R(x^i)$} of monomials are projected on the respective coefficient semisheaves 
$\phi _L(a_i)$
\resp{$\phi _R(a_i)$} or 
$\phi _L(b_i)$
\resp{$\phi _R(b_i)$} of which sections are respectively the coefficient functions $a_i$ or $b_i$.

As the semisheaf $\phi _L^{\star S_p}(M_L)$
\resp{$\phi _R^{\star S_p}(M_R)$} is defined on the algebraic semigroup
$M_L$ \resp{$M_R$} according to section 2.2, the semisheaves of monomials 
$\phi _L(x^i)$
\resp{$\phi _R(x^i)$} and the semisheaves of coefficients are also algebraic and characterized by a same set of increasing ranks, being algebraic dimensions defined from global residue degrees as developed in \cite{Pie2}.
\vskip 11pt

\subsection{Blowups of the versal deformations}

{\bbf A blowup of the versal deformation 
$D_{S_L}$
\resp{$D_{S_R}$}} can be envisaged: it {\bbf consists in the extension of the quotient algebra of the versal deformation} and corresponds to the inverse versal deformation
$D^{-1}_{S_L}$
\resp{$D^{-1}_{S_R}$}.  {\bbf It is based on the following smooth endomorphism\/}:
\begin{align*}
E_{x^i}[\phi _L(x^i)] &= \phi _L(x^i)_r\oplus\phi _L(x^i)_I\;, && 1\le i\le k-1\; ,\\
\rresp{E_{x^i}[\phi _R(x^i)] &= \phi _R(x^i)_r\oplus\phi _R(x^i)_I}\;,
\end{align*}
{\bbf based on Galois antiautomorphisms} \cite{Pie2}, \cite{Pie8}, where:
\Bi
\item $\phi _L(x^i)_r$\resp{$\phi _R(x^i)_r$} is the residual monomial semisheaf on the respective coefficient semisheaf;

\item $\phi _L(x^i)_I$\resp{$\phi _R(x^i)_I$} is the complementary monomial semisheaf disconnected from the respective coefficient  semisheaf on which it was projected.
\Ei
\vskip 11pt

Let $\Pi _{S_L}$
\resp{$\Pi _{S_R}$} denote the set
\[ \{E_{x^i}[\phi _L(x^i)]\}_{i=1}^{k-1}
\rresp{\{E_{x^i}[\phi _R(x^i)]\}_{i=1}^{k-1}}\]
of smooth endomorphisms disconnecting totally the monomial semisheaves 
$\phi _L(x^i)$
\resp{\linebreak $\phi _R(x^i)$} from the respective coefficient semisheaves.

Let $T_{V_{x_L}}=\{\dots,T_{V_{x^i}},\dots\}$
\resp{$T_{V_{x_R}}=\{\dots,T_{V_{x^i}},\dots\}$} denote the set of tangent bundles obtained by projecting all the disconnected monomial base semisheaves 
$\phi _L(x^i)_I$
\resp{$\phi _R(x^i)_I$} in the vertical tangent space.

Then, {\bbf the extension of the quotient algebra of the versal deformation of the singular semisheaf}
$\phi ^{\star S_p} _L(M_L)$
\resp{$\phi ^{\star S_p} _R(M_R)$}, having an isolated degenerate singularity on each section, {\bbf is realized by the spreading out isomorphism}
\[
(\SOT)_L=(T_{V_{x_L}}\circ \Pi _{S_L})
\rresp{(\SOT)_R=(T_{V_{x_R}}\circ \Pi _{S_R})}\;.\]
as developed in \cite{Pie2} and \cite{Pie8}.
\vskip 11pt

Let $\phi _L(x^i)$
\resp{$\phi _R(x^i)$} and
$\phi _L(x^j)$
\resp{$\phi _R(x^j)$} be two monomial base semisheaves of the base
$\phi _{S_L}$ of the versal deformation.

Then, $\phi _L(Dx_L^{i-j})$
\resp{$\phi _R(Dx_R^{i-j})$} will denote the gluing up of these two monomial semisheaves on a connected domain 
$Dx_L^{i-j}$
\resp{$Dx_R^{i-j}$}.

Let $\phi _{\SOT(1)_L}^{S_p}$
\resp{$\phi _{\SOT(1)_R}^{S_p}$} be the set
$\phi _{S_L}$
\resp{$\phi _{S_R}$} of {\bbf the base monomial semisheaves totally disconnected} from 
$\phi _L^{\star S_p}$
\resp{$\phi _R^{\star S_p}$} {\bbf by the blowup 
$D_{S_L}^{-1}$
\resp{$D_{S_R}^{-1}$} of the versal deformation} in such a way that these monomial semisheaves:
\Bean
\item {\bbf are glued together section by section\/};
\item {\bbf cover partially the residue singular semisheaf
$\phi _L^{\star S_p}$
\resp{$\phi _R^{\star S_p}$}\/}, in the sense that each section of
$\phi _L^{\star S_p}$
\resp{$\phi _R^{\star S_p}$} is totally or partially covered by the corresponding section of
$\phi _{\SOT(1)_L}^{S_p}$
\resp{$\phi _{\SOT(1)_R}^{S_p}$} obtained by gluing up the base monomials of the versal deformation.
\Ee
\vskip 11pt

Remark that the {\bbf blowup 
$D_{S_L}^{-1}$
\resp{$D_{S_R}^{-1}$}} of the versal deformation {\bbf has been envisaged as being maximal\/}, i.e. when the base monomial semisheaves are totally disconnected from the singular semisheaf
$\phi _L^{\star S_p}$
\resp{$\phi _R^{\star S_p}$}.  The intermediate cases, i.e. when the monomial semisheaves are partially disconnected from
$\phi _L^{\star S_p}$
\resp{$\phi _R^{\star S_p}$}, are studied in \cite{Pie2} and in \cite{Pie8}.
\vskip 11pt

\subsection{Sequence of blowups of versal deformations}

Referring to section 3.3, it appears that the blowup of the versal deformation of the swallowtail $F(x,a_1,a_2,a_3)=x^5+a_1x^1+a_2x^2+a_3x^3$ generates especially the singular base monomial $f(x)=x^3$.  Consequently, {\bbf the blowup of the versal deformation of the singular semisheaf}
$\phi _L^{\star S_p}$
\resp{$\phi _R^{\star S_p}$} {\bbf of which sections are affected by degenerate singularities of type swallowtail generates monomial base semisheaves
$\phi _{\SOT(1)_L}^{S_p}$
\resp{$\phi _{\SOT(1)_R}^{S_p}$} of which
$\phi _L^{\star}(x^3)$
\resp{$\phi _R^{\star}(x^3)$}\/}, case $i=3$, {\bbf is again a singular semisheaf} noted
$\phi _L^{\star S_p}(x^3)$
\resp{$\phi _R^{\star S_p}(x^3)$}.
\vskip 11pt

This singular semisheaf
$\phi _L^{\star S_p}(x^3)$
\resp{$\phi _R^{\star S_p}(x^3)$} can then be submitted to a versal deformation and a blowup of it generating the monomial base semisheaf
$\phi _{L(2)}(x)$
\resp{$\phi _{R(2)}(x)$} covering partially by patches the semisheaf
$\phi _{\SOT(1)_L}^{S_p}$
\resp{$\phi _{\SOT(1)_R}^{S_p}$}.

Thus, {\bbf in the case of a singular semisheaf
$\phi _L^{\star S_p}(M_L)$
\resp{$\phi _R^{\star S_p}(M_R)$} of swallowtail type, the two semisheaves}
\begin{align*}
\phi _{\SOT(1)_L}^{S_p} & \rresp{\phi _{\SOT(1)_R}^{S_p}} \\
\text{and} \quad
\phi _{\SOT(2)_L}^{S_p}\equiv \phi _{L(2)}(x) &
\rresp{\phi _{\SOT(2)_R}^{S_p}\equiv \phi _{R(2)}(x)}\;, \end{align*}
generated by versal deformation and blowup of this one,
{\bbf cover partially the residual singular semisheaf}
$\phi _L^{\star S_p}(M_L)$
\resp{$\phi _R^{\star S_p}(M_R)$} {\bbf leading to the embedding\/}:
\begin{align*}
\phi _L^{\star S_p}(M_L) & \subset \phi _{\SOT(1)_L}^{S_p} \subset \phi _{\SOT(2)_L}^{S_p}\\
\rresp{\phi _R^{\star S_p}(M_R) & \subset \phi _{\SOT(1)_R}^{S_p} \subset \phi _{\SOT(2)_R}^{S_p}}\;.\end{align*}
\vskip 11pt

Referring to section 3.1, the semisheaf
$\phi _L^{\star S_p}(M_L)$
\resp{$\phi _R^{\star S_p}(M_R)$} is a space  semifield.  The covering semifields
$\phi _{\SOT(1)_L}^{S_p}$
\resp{$\phi _{\SOT(1)_R}^{S_p}$} and
$\phi _{\SOT(2)_L}^{S_p}$
\resp{$\phi _{\SOT(2)_L}^{S_p}$}
are thus also of space nature. 

The corresponding semifields of ``time'' type can be generated by a composition of morphisms 
$\gamma _{r\to t}\circ E_L$
\resp{$\gamma _{r\to t}\circ E_R$} as described in proposition 2.7:
\begin{align*}
\gamma _{r\to t}\circ E_L : \quad \phi _L^{\star S_p}(M_L) &\To\phi _L^{T_p}(M_L) \oplus \phi _L^{*S_p}(M_L)\\
\rresp{\gamma _{r\to t}\circ E_R : \quad \phi _R^{\star S_p}(M_R) &\To\phi _R^{T_p}(M_R) \oplus \phi _R^{*S_p}(M_R)}\;.\end{align*}
Let us use  the more condensed notation
$\wt M_{ST_L}^{T_p-S_p}$
\resp{$\wt M_{ST_R}^{T_p-S_p}$} for
$\phi _L^{T_p}(M_L) \oplus \phi _L^{*S_p}(M_L)$
\resp{$\phi _R^{T_p}(M_R) \oplus \phi _R^{*S_p}(M_R)$}.

Similarly, time semisheaves can be generated by
$(\gamma _{r\to t}\circ E_L)$
\resp{$(\gamma _{r\to t}\circ E_R)$} morphisms from the space  semisheaves
$\phi _{\SOT(1)_L}^{S_p}$
\resp{$\phi _{\SOT(1)_R}^{S_p}$} and
$\phi _{\SOT(2)_L}^{S_p}$
\resp{$\phi _{\SOT(2)_R}^{S_p}$} leading to covering space-time semisheaves
$\phi _{\SOT(1)_L}^{T_p-S_p}$
\resp{$\phi _{\SOT(1)_R}^{T_p-S_p}$} and
$\phi _{\SOT(2)_L}^{T_p-S_p}$
\resp{$\phi _{\SOT(2)_R}^{T_p-S_p}$}.

{\bbf These two covering space-time semifields will be respectively labeled  ``$\MG$''  and  ``$\M$''\/}, as introduced in \cite{Pie8}, and will be noted according to:
\begin{align*}
\wt M_{\MG_L}^{T_p-S_p} &\equiv \phi _{\SOT(1)_L}^{T_p-S_p} \qquad 
\rresp{\wt M_{\MG_R}^{T_p-S_p} \equiv \phi _{\SOT(1)_R}^{T_p-S_p}}\phantom{\;,}&     \\
\text{and} \quad \wt M_{\M_L}^{T_p-S_p} &\equiv \phi _{\SOT(2)_L}^{T_p-S_p} \qquad
\rresp{\wt M_{\M_R}^{T_p-S_p} \equiv \phi _{\SOT(2)_R}^{T_p-S_p}} \;.&\end{align*}
Thus, {\bbf degenerate singularities of corank 1 and codimension 3 on the semisheaf}\linebreak
$\wt M_{\ST_L}^{S_p}\equiv \phi _L^{\star S_p}(M_L)$
\resp{$\wt M_{\ST_R}^{S_p}\equiv \phi _R^{\star S_p}(M_R)$} {\bbf are able to generate the three embedded semisheaves\/}:
\[
\wt M_{\ST_L}^{T_p-S_p}  \subset \wt M_{\MG_L}^{T_p-S_p} \subset  \wt M_{\M_L}^{T_p-S_p} \qquad
\rresp{\wt M_{\ST_R}^{T_p-S_p}  \subset \wt M_{\MG_R}^{T_p-S_p} \subset  \wt M_{\M_R}^{T_p-S_p}}\;.\]
\vskip 11pt

Remark that instead of considering degenerate singularities on a  semifield of ``space type''\linebreak
$\wt M_{\ST_L}^{S_p}\equiv \phi _L^{\star S_p}(M_L)$, we could have envisaged them on the corresponding  semifield of ``time type''
$\wt M_{\ST_L}^{T_p}\equiv \phi _L^{\star T_p}(M_L)$ leading similarly to the same three embedded semisheaves of type $\ST\subset \MG\subset \M$ as described explicitly in \cite{Pie8}.
\vskip 11pt

\subsection{Proposition (Langlands general global correspondences)}

{\em
Let
\[
\sigma ( W_{\wt L_{\o v}}\times W_{\wt L_{v}} )_{\ST}
= \Repsp ( ( \GL_2^{(r)} ( L_{\o v}\times L_v) )  +( \GL_2^{\perp} ( L_{\o v}\times L_v ) ))_{\ST}
\]
be the representation space of the product, right by left, of global Weil groups given by the time and space orthogonal bilinear algebraic semigroups ``$\ST$'' 
$\GL_2^{(r)}(L_{\o v}\times L_v))_{\ST}$ and $\GL_2^{\perp}(L_{\o v}\times L_v))_{ST}$.

Let 
\begin{align*}
&\sigma ( W_{\wt L_{\o v_{\SOT(1)}}}\times W_{\wt L_{v_{\SOT(1)}}})_{\MG}\\
& \qquad \qquad = \Repsp ( ( \GL_2^{(r)} ( L_{\o v_{\cov(1)}}\times L_{v_{\cov(1)}}))+(\GL_2^{\perp}(L_{\o v_{\cov(1)}}\times L_{v_{\cov(1)}})))_{\MG}\\
\rresp{&\sigma ( W_{\wt L_{\o v_{\SOT(2)}}}\times W_{\wt L_{v_{\SOT(2)}}} )_{\M}\\
& \qquad \qquad = \Repsp ( ( \GL_2^{(r)}(L_{\o v_{\cov(2)}}\times L_{v_{\cov(2)}}))+( \GL_2^{\perp}( L_{\o v_{\cov(2)}}\times L_{v_{\cov(2)}})))_{\M}}
\end{align*}
be the representation space of the product, right by left, of covering global Weil groups given by the covering time and space bilinear algebraic semigroups ``$\MG$'' \resp{``$\M$''}, where
$L_{\o v_{\cov(1)}}$ \resp{$L_{\o v_{\cov(2)}}$} is the set of completions covering $L_v$ and corresponding to the set of extensions $\wt L_{\o v_{\SOT(1)}}$ \resp{$\wt L_{\o v_{\SOT(2)}}$} of which degrees are inferior or equal to that of
$L_v$ \cite{Pie2}.

Let $(\wt M_{\ST_R}^{T_p-S_p}\otimes\wt M_{\ST_L}^{T_p-S_p})$,
$(\wt M_{\MG_R}^{T_p-S_p}\otimes\wt M_{\MG_L}^{T_p-S_p})$ and
$(\wt M_{\M_R}^{T_p-S_p}\otimes\wt M_{\M_L}^{T_p-S_p})$ denote respectively the 
``$\ST$'',
``$\MG$'' and
``$\M$'' bisemisheaves of differentiable bifunctions on the time and space bilinear algebraic semigroups
``$\ST$'',
``$\MG$'' and
``$\M$'' defined above.

{\bbf After desingularization {\em \cite{Abh}, \cite{DeJ}} of the bisemisheaves
$(\wt M_{\ST_R}^{T_p-S_p}\otimes\wt M_{\ST_L}^{T_p-S_p})$ and
$(\wt M_{\MG_R}^{T_p-S_p}\otimes\wt M_{\MG_L}^{T_p-S_p})$
and toroidal compactification of the bisemisheaves
$(\wt M_{\ST_R}^{T_p-S_p}\otimes\wt M_{\ST_L}^{T_p-S_p})$,
$(\wt M_{\MG_R}^{T_p-S_p}\otimes\wt M_{\MG_L}^{T_p-S_p})$ and
$(\wt M_{\M_R}^{T_p-S_p}\otimes\wt M_{\M_L}^{T_p-S_p})$, we get the following general global correspondences of Langlands\/}:
\begin{align*}
\LGGC_{\ST} : \quad &\sigma (W_{\wt L_{\o v}}\otimes W_{\wt L_v}) \\
&\qquad  \overset{\sim}{\To} \Pi (\GL_2^{(r)}(L_{\o v}\times L_v))+\Pi (\GL_2^\perp(L_{\o v}\times L_v))\\
\LGGC_{\MG} : \quad &\sigma (W_{\wt L_{\o v_{\SOT(1)}}}\otimes W_{\wt L_{v_{\SOT(1)}}}) \\
&\qquad  {\To} \Pi (\GL_2^{(r)}(L_{\o v_{\cov(1)}}\times L_{v_{\cov(1)}}))+\Pi (\GL_2^\perp(L_{\o v_{\cov(1)}}\times L_{v_{\cov(1)}}))\\
\LGGC_{\M} : \quad &\sigma (W_{\wt L_{\o v_{\SOT(2)}}}\otimes W_{\wt L_{v_{\SOT(2)}}}) \\
&\qquad  {\To} \Pi (\GL_2^{(r)}(L_{\o v_{\cov(2)}}\times L_{v_{\cov(2)}}))+\Pi (\GL_2^\perp(L_{\o v_{\cov(2)}}\times L_{v_{\cov(2)}}))
\end{align*}
where $\Pi (\GL_2^{(r)}(\cdots \times \cdots))$ denotes the cuspidal representation of the considered bilinear algebraic semigroup.
}
\vskip 11pt

\begin{proof}
We refer to the notations of proposition 2.8 where the Langlands general global correspondence
$\LGGC_{\ST}$ was already introduced.

And, more particularly, we have that:
\Bi
\item the representation space
$\sigma (W_{\wt L_{\o v_{\L\{\begin{subarray}{c} {\sSOT}(1)\\ \sSOT(2)\end{subarray}\R.}}}\otimes W_{\wt L_{v_{\L\{\substack{\sSOT(1)\\ \sSOT(2)}\R.}}})$
of the product, right by left, of global Weil groups decomposes into a reduced and an orthogonal part according to:
\begin{multline*}
\sigma (W_{\wt L_{\o v_{\L\{\begin{subarray}{c} {\sSOT}(1)\\ \sSOT(2)\end{subarray}\R.}}} \otimes W_{\wt L_{v_{\L\{\begin{subarray}{c} {\sSOT}(1)\\ \sSOT(2)\end{subarray}\R.}}})
\\
=
\sigma^{(r)} (W_{\wt L_{\o v_{\L\{\begin{subarray}{c} {\sSOT}(1)\\ \sSOT(2)\end{subarray}\R.}}}\otimes W_{\wt L_{v_{\L\{\begin{subarray}{c} {\sSOT}(1)\\ \sSOT(2)\end{subarray}\R.}}})\oplus
\sigma^{\perp} (W_{\wt L_{\o v_{\L\{\begin{subarray}{c} {\sSOT}(1)\\ \sSOT(2)\end{subarray}\R.}}}\otimes W_{\wt L_{v_{\L\{\begin{subarray}{c} {\sSOT}(1)\\ \sSOT(2)\end{subarray}\R.}}})\;;
\end{multline*}

\item the bisemisheaves 
$(\wt M_{\ST_R}^{T_p-S_p}\otimes\wt M_{\ST_L}^{T_p-S_p})$,
$(\wt M_{\MG_R}^{T_p-S_p}\otimes\wt M_{\MG_L}^{T_p-S_p})$ and
$(\wt M_{\M_R}^{T_p-S_p}\otimes\wt M_{\M_L}^{T_p-S_p})$ are functional representation spaces of the corresponding bilinear algebraic semigroups.
\Ei

As it results from section 3.5, the bisemisheaves
$(\wt M_{\ST_R}^{T_p-S_p}\otimes\wt M_{\ST_L}^{T_p-S_p})$ and
$(\wt M_{\MG_R}^{T_p-S_p}\otimes\wt M_{\MG_L}^{T_p-S_p})$ 
may be affected by residual singularities after the envisaged versal deformations and blowups of these.  Consequently, these bisemisheaves must be desingularized according to the classical procedure recalled in \cite{Pie2} in order to envisage a cuspidal representation of these.

Furthermore, a toroidal compactification of the bisemisheaves
$(\wt M_{\ST_R}^{T_p-S_p}\otimes\wt M_{\ST_L}^{T_p-S_p})$,
$(\wt M_{\MG_R}^{T_p-S_p}\otimes\wt M_{\MG_L}^{T_p-S_p})$ and
$(\wt M_{\M_R}^{T_p-S_p}\otimes\wt M_{\M_L}^{T_p-S_p})$ must be undertaken on these, as it was done in proposition 2.8, to reach the searched cuspidal representations.
\end{proof}
\vskip 11pt

\subsection{Proposition (The three singular and nonsingular universal $\GL_2$-structures)}

{\em 
\Bena
\item The singular bisemisheaves 
$(\wt M_{\ST_R}^{T_p-S_p}\otimes\wt M_{\ST_L}^{T_p-S_p})$ and
$(\wt M_{\MG_R}^{T_p-S_p}\otimes\wt M_{\MG_L}^{T_p-S_p})$ and the nonsingular bisemisheaf
$(\wt M_{\M_R}^{T_p-S_p}\otimes\wt M_{\M_L}^{T_p-S_p})$ are universal mathematical structures.

\item
\Be
\item The direct sum of the ``$\ST$ and ``$\MG$''  bisemisheaves
$(\wt M_{\ST_R}^{T_p-S_p}\oplus\wt M_{\MG_R}^{T_p-S_p}) \otimes
(\wt M_{\ST_L}^{T_p-S_p}\oplus\wt M_{\MG_L}^{T_p-S_p})$, including their interactions, is a two-level universal physical structure corresponding to the space-time string fields of the dark matter structure of an elementary particle.

\item The direct sum of the ``$\ST$, ``$\MG$''  and ``$\M$'' bisemisheaves
$(\wt M_{\ST_R}^{T_p-S_p}\oplus\wt M_{\MG_R}^{T_p-S_p}\oplus\wt M_{\M_R}^{T_p-S_p}) \otimes
(\wt M_{\ST_L}^{T_p-S_p}\oplus\wt M_{\MG_L}^{T_p-S_p}\oplus\wt M_{\M_L}^{T_p-S_p})$, inclusing their interactions, is a three-level universal physical structure corresponding to the space-time string fields of the visible matter structure of an elementary particle.
\Ee
\Ee}
\vskip 11pt

\begin{proof}
\Bena
\item Taking into account the desingularization of the ``$\ST$ and ``$\MG$''  bisemisheaves
$(\wt M_{\ST_R}^{T_p-S_p}\otimes\wt M_{\ST_L}^{T_p-S_p})$ and
$(\wt M_{\MG_R}^{T_p-S_p}\otimes\wt M_{\MG_L}^{T_p-S_p})$, they are, as the nonsingular 
``$\M$  bisemisheaf
$(\wt M_{\M_R}^{T_p-S_p}\otimes\wt M_{\M_L}^{T_p-S_p})$, universal mathematical structures as developed in proposition 2.9.

Indeed, these three (desingularized) bisemisheaves, generated from the products, right by left, of global Weil groups, are in one-to-one correspondence with the holomorphic and cuspidal representations of the associated $\GL_2$-bilinear algebraic semigroups by means of the Langlands general global correspondences
$\LGGC_\ST$,
$\LGGC_\MG$ and
$\LGGC_\M$ of proposition 3.6.

In fact, the sums of the bisections of these bisemisheaves on the toroidal compactified conjugacy class representatives of the considered bilinear algebraic semigroups cover the products, right by left, of cusp forms as developed in \cite{Pie5}: in this sense, these bisemisheaves are cuspidal representations.

\item
\Be
\item Let the ``$\MG$ bisemisheaf
$(\wt M_{\MG_R}^{T_p-S_p}\otimes\wt M_{\MG_L}^{T_p-S_p})$ be generated from the ``$\ST$  bisemisheaf
$(\wt M_{\ST_R}^{T_p-S_p}\otimes\wt M_{\ST_L}^{T_p-S_p})$ by versal deformations and blowups of degenerate singularities of corank $1$ and codimension $1$ and $2$ or of degenerate singularities of corank $2$ and codimension $3$ on it.

Then,
\begin{align*}
\DMS_R(F_R) \otimes
\DMS_L(F_L) 
&= (\wt M_{\ST_R}^{T_p-S_p}\oplus\wt M_{\MG_R}^{T_p-S_p})\otimes
(\wt M_{\ST_L}^{T_p-S_p}\oplus\wt M_{\MG_L}^{T_p-S_p})\\
&= [(\wt M_{\ST_R}^{T_p-S_p}\otimes\wt M_{\ST_L}^{T_p-S_p})\oplus
(\wt M_{\MG_R}^{T_p-S_p}\otimes\wt M_{\MG_L}^{T_p-S_p})]\\
& \qquad \oplus
[(\wt M_{\ST_R}^{T_p-S_p}\otimes\wt M_{\MG_L}^{T_p-S_p})\oplus
(\wt M_{\MG_R}^{T_p-S_p}\otimes\wt M_{\ST_L}^{T_p-S_p})]_{\int}
 \end{align*}
corresponds to the dark matter (string) fields of an elementary bisemifermion $(F_R\otimes F_L)$ where the bisemisheaves inside $[\cdot]_{\int}$ are interaction semifields between mixed right and left ``$\ST$ and ``$\MG$'' semisheaves.

The bisemisheaf $(\wt M_{\MG_R}^{T_p-S_p}\otimes\wt M_{\MG_L}^{T_p-S_p})$ corresponds to space-time
``$\MG$'' internal fields of an elementary (bisemi)fermion for the same reasons for which the
``$\ST$'' bisemisheaves $(\wt M_{\ST_R}^{T_p-S_p}\otimes\wt M_{\ST_L}^{T_p-S_p})$ are internal fields as developed in proposition 2.9.

This ``$\MG$'' bisemisheaf is a ``mass'' structure field because it is of contracting nature limiting the expansion of the internal vacuum dark energy ``$\ST$'' substructure
$(\wt M_{\ST_R}^{T_p-S_p}\otimes\wt M_{\ST_L}^{T_p-S_p})$.

As this ``$\MG$'' bisemisheaf is a mass field localized in the own vacuum of an elementary (bisemi)fermion, it must correspond to a dark mass substructure covering the dark energy substructure
$(\wt M_{\ST_R}^{T_p-S_p}\otimes\wt M_{\ST_L}^{T_p-S_p})$ of an elementary bisemifermion.

Thus, the direct sum of the ``$\ST$ and ``$\MG$''  bisemisheaves (or fields) corresponds to the dark matter structure of an elementary (bisemi)fermion.

\item Let the ``$\MG$ and ``$\M$''  bisemisheaves
$(\wt M_{\MG_R}^{T_p-S_p}\otimes\wt M_{\MG_L}^{T_p-S_p})$ and
$(\wt M_{\M_R}^{T_p-S_p}\otimes\wt M_{\M_L}^{T_p-S_p})$ be generated from the  ``$\ST$ bisemisheaf
$(\wt M_{\ST_R}^{T_p-S_p}\otimes\wt M_{\ST_L}^{T_p-S_p})$ by versal deformations and blowups of degenerate singularities of corank $1$ and codimension $3$ on it.

Then,
\begin{multline*}
\VMS_R(F_R)\otimes \VMS_L(F_L)\\
\begin{aligned}
&= (\wt M_{\ST_R}^{T_p-S_p}\oplus\wt M_{\MG_R}^{T_p-S_p}\oplus\wt M_{\M_R}^{T_p-S_p})\otimes
 (\wt M_{\ST_L}^{T_p-S_p}\oplus\wt M_{\MG_L}^{T_p-S_p}\oplus\wt M_{\M_L}^{T_p-S_p})\\
 &= 
 (\wt M_{\ST_R}^{T_p-S_p}\otimes\wt M_{\ST_L}^{T_p-S_p})\oplus
 (\wt M_{\MG_R}^{T_p-S_p}\otimes\wt M_{\MG_L}^{T_p-S_p})\oplus
 (\wt M_{\M_R}^{T_p-S_p}\otimes\wt M_{\M_L}^{T_p-S_p})\\
& \qquad \oplus [
 (\wt M_{\ST_R}^{T_p-S_p}\otimes\wt M_{\MG_R}^{T_p-S_p})\oplus \dots \oplus
 (\wt M_{\M_R}^{T_p-S_p}\otimes\wt M_{\MG_L}^{T_p-S_p})]
_{\int}
 \end{aligned}\end{multline*}
 corresponds to the visible matter stringfields of an elementary bisemifermion $(F_R\otimes F_L)$ where:
 \Bi
 \item $(\wt M_{\M_R}^{T_p-S_p}\otimes\wt M_{\M_L}^{T_p-S_p})$ denotes the visible mass covering substructure ``$\M$'' of this bisemifermion;
 \item the bisemisheaves inside $[\cdot]_{\int}$ are interaction fields between the three different right and left semifields  ``$\ST$, ``$\MG$''  and ``$\M$''.
 \Ei
 
 This visible mass substructure
 $(\wt M_{\M_R}^{T_p-S_p}\otimes\wt M_{\M_L}^{T_p-S_p})$ is a ``mass'' structure field because it is of contracting nature stabilizing the internal vacuum substructure, i.e. the dark matter structure
 $\DMS_R(F_R)\otimes \DMS_L(F_L)$.
 
 This mass substructure ``$\M$'' is visible because the frequency (of rotation) of its (bi)sections is inferior to that of the dark energy and dark mass substructures taking into account that the number of (bi)quanta on the (bi)sections of the bisemisheaf ``$\M$'' is inferior to the number of (bi)quanta on the respective (bi)sections of the ``$\MG$'' and ``$\ST$'' semisheaves.
 
 \item The visible matter structure of elementary particles is composed of the three-embedded shells
 $\ST \subset \MG \subset \M$ in such a way that the external shell
 ``$\M$'' is generated from the middle ground  shell ``$\MG$'' on the basis of degenerate singularities of corank $1$ and codimension $3$ on the vacuum ``$\ST$'' shell, when ``vacuum'' dark matter structure is composed of two embedded shell $\ST \subset \MG$ of which ``$\MG$'' shell is generated from degenerate singularities of corank $1$ or $2$ and codimension $<3$ on the ``$\ST$'' shell.
 
 The generation of the visible mass shell ``$\M$'' needs thus much more energy than the generation of the dark mass shell ``$\MG$'': this explains the preponderance of dark matter over ordinary visible matter in our universe.\qedhere
 \Ee
 \Ee
 \end{proof}
 \vskip 11pt
 
\subsection{Proposition (Sets of Langlands general global correspondences)}

{\em
According to the existence or absence of degenerate singularities on the basic space-time bisemisheaf
$(\wt M_{\ST_R}^{T_p-S_p}\otimes\wt M_{\ST_L}^{T_p-S_p})$, there are three sets of Langlands general global correspondences ($\LGGC$).
\Bean
\item If there are no degenerate singularities on
$(\wt M_{\ST_R}^{T_p-S_p}\otimes\wt M_{\ST_L}^{T_p-S_p})$, we get the one level $\LGGC_\ST$:
\[ \LGGC_\ST : \qquad
\sigma (W_{\wt L_{\o v}}\times W_{\wt L_v})
\To \Pi (\GL_2^{(r)}(L_{\o v}\times L_v))+ \Pi (\GL_2^{\perp}(L_{\o v}\times L_v))\;.\]

\item If there are degenerate singularities of corank $1$ and codimension inferior to $3$ or degenerate singularities of corank $2$ and codimension inferior or equal to $3$ on
$(\wt M_{\ST_R}^{T_p-S_p}\otimes\wt M_{\ST_L}^{T_p-S_p})$, we get the two-level $\LGGC_{\ST-\MG}$:
\begin{multline*}
\LGGC_{\ST-\MG}: \qquad \sigma (W_{\wt L_{\o v}}\times W_{\wt L_v})+
\sigma (W_{\wt L_{\o v_{\SOT(1)}}}\times W_{\wt L_{v_{\SOT(1)}}})\\
\To \begin{aligned}[t]
&\Pi (\GL_2^{(r)}(L_{\o v}\times L_v))+\Pi (\GL_2^{\perp}(L_{\o v}\times L_v))\\
& \quad + \Pi (\GL_2^{(r)}(L_{\o v_{\cov(1)}}\times L_{v_{\cov(1)}}))+\Pi (\GL_2^{\perp}(L_{\o v_{\cov(1)}}\times L_{v_{\cov(1)}}))\;.
\end{aligned}\end{multline*}

\item If there are degenerate singularities of corank $1$ and codimension $3$ on
$(\wt M_{\ST_R}^{T_p-S_p}\otimes\wt M_{\ST_L}^{T_p-S_p})$, we get the three level $\LGGC_{\ST-\MG-\M}$:
\begin{multline*}
\LGGC_{\ST-\MG-\M}: \quad \sigma (W_{\wt L_{\o v}}\times W_{\wt L_v})+
\sigma (W_{\wt L_{\o v_{\SOT(1)}}}\times W_{\wt L_{v_{\SOT(1)}}})+
\sigma (W_{\wt L_{\o v_{\SOT(2)}}}\times W_{\wt L_{v_{\SOT(2)}}})\\
\To \begin{aligned}[t]
&\Pi (\GL_2^{(r)}(L_{\o v}\times L_v))+\Pi (\GL_2^{\perp}(L_{\o v}\times L_v))\\
& \quad + \Pi (\GL_2^{(r)}(L_{\o v_{\cov(1)}}\times L_{v_{\cov(1)}}))+\Pi (\GL_2^{\perp}(L_{\o v_{\cov(1)}}\times L_{v_{\cov(1)}}))\\
& \quad + \Pi (\GL_2^{(r)}(L_{\o v_{\cov(2)}}\times L_{v_{\cov(2)}}))+\Pi (\GL_2^{\perp}(L_{\o v_{\cov(2)}}\times L_{v_{\cov(2)}}))\;.
\end{aligned}\end{multline*}
\Ee}
\vskip 11pt

\begin{proof}
These three sets of Langlands general global correspondences directly result from propositions 3.6 and 3.7 in such a way that the one-level correspondence $\LGGC_\ST$ is related to a nonsingular universal $\GL(2)$-structure while the two- and three-level correspondences 
$\LGGC_{\ST-\MG}$ and
$\LGGC_{\ST-\MG-\M}$ are related to singular universal $\GL_2$-structures.
\end{proof}

\vfill
C. Pierre\\
Universit\'e de Louvain\\
Chemin du cyclotron, 2\\
B--1348 Louvain-la-Neuve, Belgium\\
pierre.math.be@gmail.com

\end{document}